\documentclass[dvips,twoside]{article}
\usepackage[a5paper,layoutwidth=148.5mm,layoutheight=7.77817in,
outer=10mm,inner=38pt,tmargin=10mm,bmargin=10mm,footskip=15pt,
]{geometry}

\usepackage[utf8]{inputenc}  
\usepackage[T1]{fontenc}
\usepackage{lmodern}
\usepackage{fix-cm}
\usepackage
{microtype}
\usepackage{amsmath}
\usepackage[square]{natbib}
\setcitestyle{square,aysep={}}
\usepackage[french]{babel}
\usepackage[autolanguage,np]{numprint}
\usepackage{graphicx}        
\usepackage[all]{xy} 

\newcommand\mapright[1]{\mathop{\longrightarrow}\limits^{#1}}

\usepackage{authblk}

\makeatletter
\renewcommand\@spart[1]{%
  {\parindent \z@ \raggedright
    \interlinepenalty \@M
    \normalfont
    \huge \bfseries\raggedright #1\par}%
  \nobreak
  \vskip 3ex
  \@afterheading}
\renewcommand\section{\@startsection {section}{1}{\z@}%
  {-3.5ex \@plus -1ex \@minus -.2ex}%
  {2.3ex \@plus.2ex}%
  {\normalfont\Large\bfseries\raggedright}}
\renewcommand\subsection{\@startsection{subsection}{2}{\z@}%
  {-3.25ex\@plus -1ex \@minus -.2ex}%
  {1.5ex \@plus .2ex}%
  {\normalfont\large\bfseries\raggedright}}
\renewcommand\subsubsection{\@startsection{subsubsection}{3}{\z@}%
  {-3.25ex\@plus -1ex \@minus -.2ex}%
  {1.5ex \@plus .2ex}%
  {\normalfont\normalsize\bfseries\raggedright}}
\makeatother

\usepackage
{hyperref}
\hypersetup{bookmarksopen=false,breaklinks=true,
  plainpages=false,
  colorlinks=true,linkcolor=blue,
  citecolor=red,urlcolor=red,hypertexnames=false}                      

\usepackage{breakurl}

\usepackage{doi}

\begin{document}

\DeclareFontShape{T1}{lmr}{bx}{sc} { <-> ssub * cmr/bx/sc }{}

\title{
  Enquête sur les modes d'existence\\des êtres
  mathématiques (version augmentée)}

\author{Guy Wallet
  \and
  Stefan Neuwirth}

\date{
}

\maketitle

\begin{abstract}
  L'objet de cet essai est de s'interroger sur la manière dont les
  entités mathématiques pourraient être accueillies dans
  l'architecture des modes d'existence proposée par Bruno Latour dans
  le cadre de son ontologie renouvelée et pluraliste du monde moderne
  \citep{Latour2012,EME}. Après avoir décrit les termes du problème,
  les travaux de Reviel Netz sur l'émergence des mathématiques
  grecques \citep{Netz1999} et ceux de Charles Sanders Peirce sur la
  dimension diagrammatique de l'activité mathématique
  \citep{Peirce33-58,Peirce76} sont présentés, ainsi que leurs
  conséquences relativement au thème du présent essai. Sa partie
  centrale développe longuement une conception empirique des
  mathématiques qui joue un rôle essentiel dans la suite. Cette
  analyse est basée sur la notion d'expérience chère à William James
  \citep{James2007}; elle est aussi inspirée par certains aspects de la
  philosophie de Per Martin-Löf \citep{MartinLof87}. Elle permet de
  penser la solide certitude dont la démonstration dote les résultats
  mathématiques tout en invalidant l'interprétation de cette
  certitude comme la marque d'un accès direct à une vérité absolue et
  transcendante. En s'appuyant sur cette analyse, la suite de ce
  travail définit une forme de quasi-mode d'existence propre aux êtres
  mathématiques qui respecte les traits principaux d'un mode
  d'existence selon l'ontologie latourienne. En conclusion, des
  éléments de discussion sont apportés quant à la manière dont ce
  quasi-mode pourrait être placé dans cette ontologie, notamment par
  rapport au mode de la référence objective qui prévaut dans de
  nombreuses autres sciences.
\end{abstract}

\section{Introduction}

Désignée par l'abréviation \textsc{eme} dans la suite de ce texte,
l'\emph{Enquête sur les modes d'existence} est d'abord un livre de
Bruno Latour, \citep{Latour2012}. La version numérique de cet ouvrage
est étoffée par l'enquête collective \citep{EME} à laquelle les
internautes ont la possibilité de participer.  Le livre et sa version
numérique augmentée visent à renouveler notre vision des modes
d'existence et des régimes de vérité associés. Il s'agit d'un projet
de grande ampleur qui se propose de revisiter le cœur de notre
modernité --~les sciences, les techniques, le droit, la religion, la
politique, l'économie, etc.~-- et d'élaborer une description
explicite d'une multiplicité de formes d'existence, chacune dotée de
ses propres conditions de félicité. Ce projet d'une \emph{ontologie
  pluraliste} se place en opposition avec le point de vue qui
accompagne la modernité selon lequel il n'y a qu'une seule forme
d'existence, à savoir celle qui relève de la Science et plus largement
de la Raison.  Cette opposition n'a nullement pour objet de critiquer
l'institution scientifique, ni de nier son importance; le propos de ce
projet est plutôt de mettre la science à sa juste place en restituant
finement la riche spécificité de son fonctionnement tout en concevant
la possibilité d'autres institutions et d'autres formes d'existence
elles aussi significatives de la modernité\footnote{L'un des motifs de
  cette refonte ontologique est de fournir un nouveau point de vue
  approprié pour penser et réagir face à la menace d'une catastrophe
  écologique et climatique de plus en plus pressante. Malgré son
  immense intérêt, cet aspect ne sera pas évoqué dans la suite du fait
  qu'il n'interfère pas immédiatement avec le sujet traité dans le
  présent article.}.

De fait, un moment essentiel dans la genèse de l'\textsc{eme} est la
mise en évidence du mode d'existence prépondérant dans la démarche
scientifique pour produire une connaissance objective, avant tout dans
les sciences expérimentales et de terrain\footnote{Dans un premier
  temps, le premier auteur du présent texte interprétait ce mode comme
  spécifique à la démarche scientifique expérimentale et de terrain,
  et caractérisant cette dernière. Après discussion avec Isabelle
  Stengers et Bruno Latour, il semble que cette interprétation est un
  peu rapide. D'une part, on peut observer la mise en œuvre de ce mode
  en dehors du cadre strictement scientifique, comme le montre la
  place occupée par l'exemple d'une excursion sur le mont Aiguille
  dans l'\textsc{eme}. D'autre part, il n'est pas absolument clair que
  ce même mode recouvre tout ce qui se fait dans les sciences
  expérimentales et de terrain.}, ce qui laisse provisoirement à
l'écart les mathématiques. Il s'agit du \emph{mode de la référence}
--~référence objective~-- noté [\textsc{ref}] et caractérisé par la
construction de \emph{chaines de référence}. Une description un peu
plus précise de ce mode et de la notion de chaine de référence sera
donnée vers la fin de ce travail (voir la partie~\ref{PrecisREF}) au
moyen de la présentation d'un exemple caractéristique. Il suffit pour
l'instant de savoir qu'une chaine de référence est une succession (un
enchainement) de médiations (d'inscriptions scripto-visuelles), chaque
inscription\footnote{La description dans l'\textsc{eme} du mode de la
  référence fait grand usage du terme d'\emph{inscription} ou de
  manière équivalente de celui de \emph{forme}, voire de celui
  d'\emph{idéographie}, pour désigner un dispositif concret qui permet
  d'opérer une transition entre une face plus matérielle et une face
  qui relève plus de l'écriture ou du calcul. Pour citer
  l'\textsc{eme}: «par exemple une vitrine d'exposition sur laquelle
  on pose un spécimen est une inscription au même titre que l'écran
  d'un scanner médical ou un tableau excel» \cite[«Livre», colonne
  «Vocabulaire», entrée «Inscription»]{EME}.}  étant produite à partir
de la précédente. La chaine s'arrête lorsque la dernière inscription
obtenue est suffisamment claire et explicite pour éclairer le problème
étudié. Finalement, quelque chose existe objectivement lorsque les
scientifiques concernés ont pu construire une chaine de référence qui
montre cette chose.

Puisque les mathématiques sont incontestablement une composante du
champ scientifique, se pose le problème de savoir si elles relèvent de
ce mode d'existence [\textsc{ref}]. Cette question n'a pas de réponse
évidente et elle est loin d'être tranchée ni même franchement abordée
dans l'état actuel de l'\textsc{eme}. Certes, la question de
l'ontologie des mathématiques est source d'un questionnement
philosophique récurrent: \emph{de facto}, le statut des êtres
mathématiques a toujours été l'objet d'interrogations et il n'est pas
étonnant de retrouver un écho de cette pierre d'achoppement dans
l'\textsc{eme}.  En première analyse, on peut préciser quelques
aspects de cette difficulté.
\begin{itemize}
\item Ranger les mathématiques dans le mode de la référence objective
  [\textsc{ref}] tel que ce dernier est présenté dans l'\textsc{eme}
  nécessite de montrer comment l'activité mathématique relève elle
  aussi des chaines de référence. Pour cela, et parce que la
  caractéristique essentielle de la pratique des mathématiques est
  l'élaboration de démonstrations, il semble nécessaire d'établir un
  pont entre la notion de chaine de référence et celle de
  démonstration. De prime abord, cela n'est pas évident et d'ailleurs,
  ce point ne semble pas évoqué dans l'\textsc{eme}.
\item Cette difficulté est certainement liée au fait que les chaines
  de référence ont été mises en évidence dans la pratique des sciences
  possédant explicitement une base expérimentale ou de terrain,
  pratique qui se déploie dans des laboratoires qui constituent le
  cadre indispensable à la réussite de ces recherches: le laboratoire
  offre des conditions matérielles, financières, politiques et
  humaines qui constituent une dimension importante de l'activité
  scientifique. L'analogue semble nettement plus difficile à concevoir
  dans le cas d'une science formelle comme les mathématiques pour
  laquelle la base expérimentale est absente en première
  approximation, et la vie de laboratoire nettement moins
  spectaculaire. Néanmoins, les dernières décennies ont vu notablement
  s'accroitre l'importance des laboratoires comme structures de base
  de la recherche en mathématiques.
\item Outre la difficulté à pénétrer au cœur de la pratique
  mathématique, l'obstacle principal à la clarification de ce problème
  est certainement la captation philosophique, vieille comme les
  mathématiques et la philosophie grecques, qui fait de la
  démonstration mathématique la forme paradigmatique de la démarche
  rationnelle parfaite, capable de mener directement et sans perte à
  la vérité pensée comme absolue et indépendante des acteurs
  humains. Cette captation est finement présentée et mainte fois
  évoquée dans l'\textsc{eme} et aussi dans l'ouvrage \citep[p.~290]{Netz1999}
  qui sera mentionné plus loin dans ce texte (voir aussi \citep[p.~189-191]{Latour2009}). Il apparait que
  cette idée d'une méthode permettant l'accès direct, plein et entier
  à la vérité, et la force qu'elle exerce implicitement dans notre
  culture occidentale, obscurcit particulièrement la réflexion sur la
  question posée.  Ce point de vue, nommé formalisme dans
  l'\textsc{eme}, a joué un rôle décisif dans l'histoire culturelle de
  la modernité. Pour reprendre les termes de Netz dans
  \citep{Netz1999}, «cette vision a hanté la culture occidentale» sur
  le long terme. D'une certaine manière, l'\textsc{eme} apparait comme
  une critique frontale de cette conception. Cependant, force est de
  reconnaitre que, même une fois dévoilé et critiqué, cet arrière-fond
  philosophique n'en continue pas moins à troubler notre réflexion et
  à rendre délicate l'identification d'un mode d'existence adéquat
  pour les mathématiques.
\end{itemize}

Et pourtant, Bruno Latour a depuis longtemps émis l'hypothèse que le
cas des sciences formelles et abstraites comme les mathématiques
devrait pouvoir se traiter de manière peut-être encore plus évidente
que celui des sciences de laboratoire ou de terrain (voir
\citep{Latour1985}). Il suppose que c'est essentiellement le poids de
nos préjugés philosophiques, particulièrement prégnants en ce qui
concerne les mathématiques, qui nous aveugle sur ce point. Le point de
vue développé dans la suite de la présente analyse pourrait en partie
valider la pertinence de cette intuition.

Notre ambition est de déterminer un mode d'existence pour les
mathématiques qui soit à la fois conforme à l'esprit de l'\textsc{eme}
et compatible avec l'expérience pratique des mathématiciens et des
utilisateurs de mathématiques (dans la suite de ce travail, le terme
de mathématicien désigne une personne ayant acquis la compétence lui
permettant d'étudier et de produire lui-même des développements
mathématiques). Le titre du présent texte évoque les modes d'existence
au pluriel parce que notre enquête identifie un autre mode, qui fait
intervenir le calcul mécanique de l'ordinateur comme étape essentielle
d'un raisonnement (comme «hiatus» dans la «trajectoire» de l'existence
de l'être mathématique, voir la partie~\ref{PrecisREF} pour l'emploi
précis de ces deux mots); il utilise l'article déterminé, «les modes»,
pour afficher notre ambition d'une enquête ouverte et la plus large
possible.

L'idée principale est de s'attaquer d'abord à l'arrière-fond
philosophique évoqué précédemment. Pour cela, une critique de la
conception formaliste des mathématiques est menée; elle permet en
contrepoint de proposer une conception empirique de la notion de
démonstration mathématique. Ce travail, peut-être le plus conséquent
de l'article, est basé sur la notion d'expérience chère à William James
\citep{James2007}, mais elle est aussi inspirée par certains aspects de
la philosophie du mathématicien, logicien et philosophe suédois Per
Martin-Löf \citep{MartinLof87}, sans oublier de précieuses réflexions
de LudwigWittgenstein \citep{Wittgenstein2006} et John Dewey \citep{Dewey2014} d'intérêt
pour ce sujet.  Ce point de vue permet de comprendre que tout résultat
mathématique est doté de la certitude résultant de sa démonstration,
tout en étant affecté d'une forme de fragilité irréductible héritée
du fait que la démonstration est d'abord une expérience. Cela permet
d'invalider l'interprétation de cette certitude comme la marque d'un
accès direct à une vérité absolue et transcendante.
 
Cette étude préparatoire mais essentielle étant menée, il est possible
d'en déduire, somme toute assez simplement, un mode d'existence --~par
prudence, il est préférable d'évoquer plutôt un quasi-mode
d'existence~-- propre aux êtres mathématiques. Ce quasi-mode respecte
pour l'essentiel les caractéristiques d'un mode d'existence selon
l'\textsc{eme}. De plus, il présente à la fois des analogies
significatives et éventuellement quelques différences notables avec le
mode d'existence de la référence objective.  La question de
l'articulation de ce quasi-mode avec celui de la référence objective
peut alors commencer à être discutée.

La présente étude est une version augmentée de
l'article \citep{walletneuwirth19}.

\section{L'apport de R. Netz et de C.~S. Peirce: la place des
  diagrammes en mathématiques} \label{NetzPeirce}

\subsection{L'émergence des mathématiques grecques selon
  Netz} \label{emerge}

Dans \citep{Latour2009}, Bruno Latour présente un ouvrage qui
constitue selon lui la première avancée décisive dans la direction
d'une «approche non formaliste du formalisme». Il s'agit du livre de
Netz publié en 1999 sous le titre \emph{The shaping of deduction in
  Greek mathematics: a study in cognitive
  history}, \citep{Netz1999}.  Cette étude présente une
analyse de l'émergence des mathématiques grecques du point de vue de
l'histoire cognitive et plus particulièrement de la constitution de la
déduction. Ce travail d'historien s'appuie sur l'étude minutieuse des
textes mathématiques grecs encore accessibles, parfois sous la forme
de fragments, écrits en gros pendant un millénaire débutant au
5\ieme~siècle av.\ J.-C\@.  L'analyse de Netz se place en opposition
explicite avec la démarche épistémologique classique, telle que celle
de Thomas Kuhn \citep{Kuhn2008}, qui s'intéresse essentiellement aux
systèmes de croyances partagées --~les paradigmes~-- sur lesquels se
fonderait le développement des connaissances scientifiques, ainsi
qu'aux grands basculements de ces paradigmes lors des révolutions
scientifiques. Netz met en doute que le travail scientifique se fonde
réellement sur ces grands paradigmes qui apparaissent plutôt après
coup dans les explications des épistémologues. Au contraire, ce sont
pour lui les pratiques partagées qui constituent la base efficiente du
développement concret des sciences. Justement, il se propose de mettre
en lumière les pratiques partagées qui indiquent le mieux la force des
ressources cognitives spécifiques aux mathématiques grecques.

\begin{figure}[h]
  \centering \includegraphics[width=300pt]{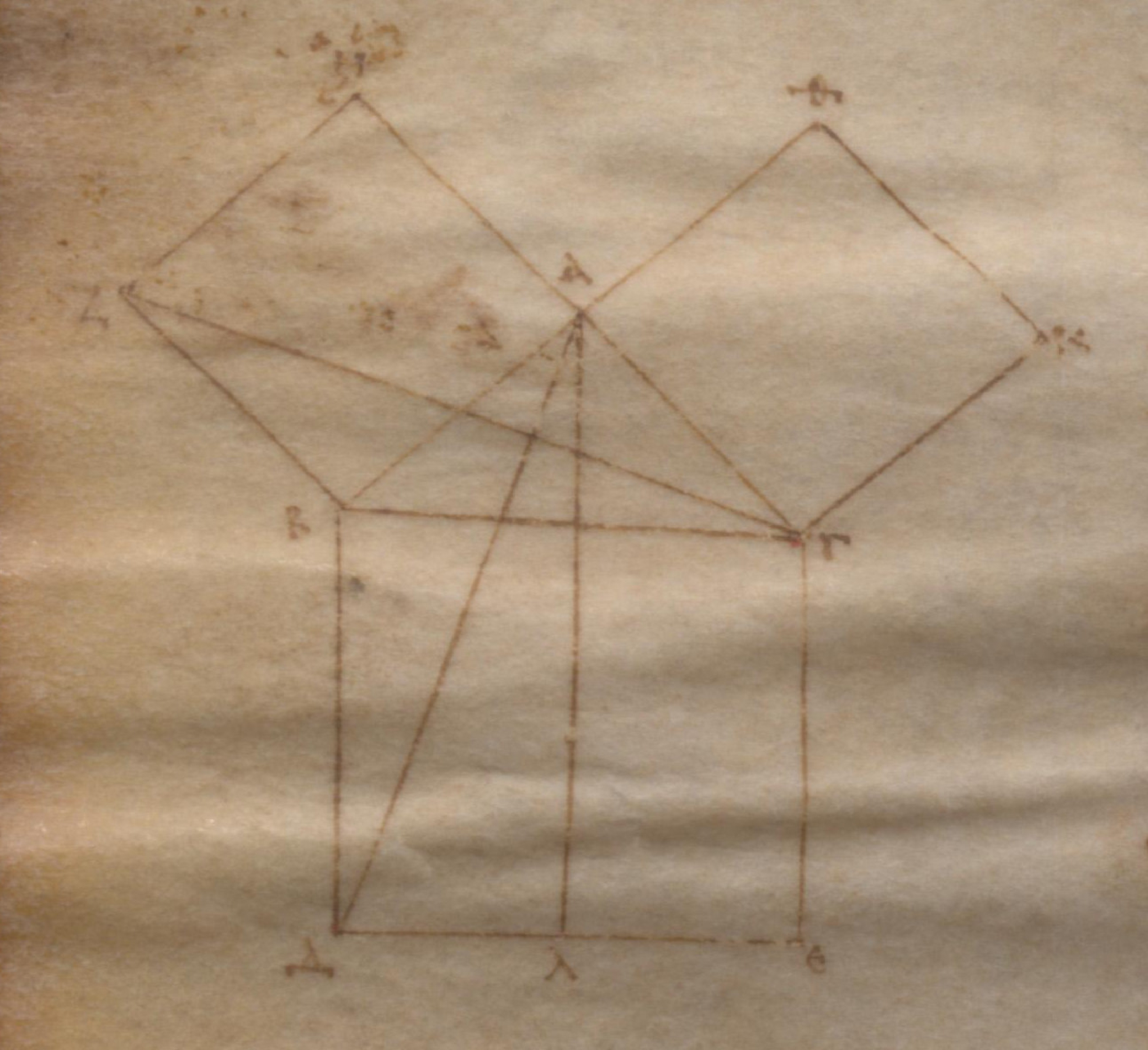}
  \caption{Le diagramme avec lettres qui accompagne le théorème de
    Pythagore dans les \emph{Éléments} d'Euclide sur le folio~31 du
    manuscrit MS D'Orville~301 de l'an~888, reproduit avec l'aimable
    autorisation de l'Institut de mathématiques Clay.}\label{fig:1}
\end{figure}

La première pratique abordée par Netz est celle de l'introduction du
diagramme avec lettres qui est l'outil le plus emblématique des
mathématiques grec\-ques: voir la figure~\ref{fig:1}.
Il s'agit à coup sûr d'un apport de première
importance qui a été, et continue à être, universellement utilisé dans
l'ensemble des sciences et des techniques. Dans le cadre des
mathématiques grecques, le diagramme avec lettres est un dispositif
qui, en relation d'interdépendance avec le texte pur, offre une vue
synoptique du cas considéré et joue le rôle d'univers de référence du
texte. Cet univers est marqué par une double finitude: les objets qui
sont marqués sont bornés et en nombre fini, ce qui en fait un bon
support pour la pensée mathématique. De plus, la nature des objets
représentés sur le diagramme n'est nulle part précisée. Non seulement
les mathématiciens grecs n'apportent aucun élément de réponse à la
question de savoir ce que pourraient bien être en vérité les entités
qui peuplent leurs raisonnements mais de plus ils ne posent jamais la
question de la réalité ultime des objets mathématiques\footnote{Un
  symptôme \emph{a contrario} en sont les difficultés rencontrées par
  un philologue et historien des mathématiques comme Árpád Szabó
  \citep{sz00} pour expliquer cette absence de questionnement.}. D'une
certaine manière, l'ambigüité ontologique semble être constitutive de
cette discipline. Ce parti pris a pour effet bénéfique d'écarter les
questions qui fâchent, celles pour lesquelles un consensus est
quasiment impossible, particulièrement dans la tradition de la
polémique propre à la démocratie grecque.

La seconde pratique des mathématiques grecques relevée par Netz
est l'utilisation d'un langage fortement spécifique. Sa première propriété est
la petite taille de son lexique: 143 mots utilisés
répétitivement sont responsables de 95\,\% du corpus
des œuvres d'Archimède, \emph{L'Arénaire} mise à part. On a déjà
vu que le diagramme avec lettres permet d'extraire de l'infinité des
entités géométriques un univers de référence fini de petite taille
suffisant et commode pour le travail de démonstration. Le lexique
mathématique fait la même chose par rapport à la vaste étendue du
lexique de la langue grecque. De plus, Netz note que ce langage
possède une structure remarquable: pour une part importante, on
l'obtient par une combinaison d'un petit nombre de briques qu'il nomme
\emph{formules}. Ce sont de petites phrases ou groupes de mots, se
répétant avec fort peu de variabilité. La présence massive de ces
formules indique qu'elles ont un rôle cognitif important dans la
structuration et la communication des
argumentaires\footnote{Cependant, la manière dont ces formules se
  combinent est complexe. Chaque démonstration met en place un réseau
  d'anaphores, doit gérer la généralité de la proposition et les
  modalités des énoncés, et Fabio Acerbi parle très à propos de
  \emph{La syntaxe logique des mathématiques grecques}
  \citep{acerbi12}.}.

En conclusion, Netz met en évidence les deux caractéristiques fortes
des mathématiques grecques que sont d'une part l'introduction du
diagramme avec lettres, qui évite tout questionnement ontologique sur
la nature ultime des objets traités, et d'autre part l'utilisation
d'un langage au lexique réduit, qui se développe implicitement selon
une «structure en formules». Selon Netz, ces deux pratiques
sont les sources principales de l'émergence et de l'efficience du mode
déductif dans les mathématiques grecques, et il en analyse en détail
la genèse.

\subsection{Le diagramme au cœur des mathématiques selon
  C.~S. Peirce} \label{Peirce}

Considéré comme le fondateur du courant pragmatiste en philosophie,
Charles Sanders Peirce est un philosophe américain. Avec Ferdinand de
Saussure, il est l'un des pères de la sémiotique (ou sémiologie),
science des signes (verbaux ou non verbaux) et du sens dont ils sont
porteurs. Peirce a été l'objet d'une reconnaissance tardive longtemps
après sa mort lorsque des chercheurs se sont mis à exploiter les
environ \nombre{12000} pages de ses écrits déposées à l'université de
Harvard. Depuis lors, ses idées suscitent un grand intérêt et de
nombreux travaux au niveau international. Dans la présente partie
consacrée à la place que Peirce accorde aux diagrammes dans les
mathématiques, les auteurs se sont principalement appuyés sur les
travaux de Christiane Chauviré \citep{Chauvire1987,Chauvire2008} pour
appréhender l'apport de ce philosophe à cet égard.
 
Le point de vue de Peirce sur les mathématiques s'appuie sur les bases
de sa sémiotique, à savoir la notion de signe et son lien avec la
pensée. En toute généralité, un signe désigne quelque chose qui,
quelle qu'en soit la raison, est interprétée comme représentant une
autre chose. Contrairement à d'autres, Peirce n'impose pas au signe
d'être simple ou minimal: les signes qu'il considère sont souvent des
agencements complexes. Par exemple, on verra bientôt qu'un diagramme
avec lettres de la géométrie grecque, aussi complexe soit-il, est un
signe au sens de Peirce au même titre qu'une simple lettre nommant un
objet géométrique dans le même diagramme.
 
Ce philosophe considère que la pensée nécessite l'utilisation des
signes au sens fort suivant: la fonction des signes n'est pas de
représenter une pensée déjà produite par un processus se déroulant
ailleurs, par exemple dans l'esprit d'un sujet humain; au contraire,
\emph{les signes constituent le substrat indispensable à la
  constitution et au développement de la pensée et, ultimement, l'acte
  de penser consiste à manipuler des signes.} De plus, Peirce
classifie les signes en trois catégories, à savoir le symbole,
l'indice et l'icône, en fonction du type de rapport qu'un
signe entretient avec son objet.  Le \emph{symbole} est un signe de nature
arbitraire dont la seule relation avec la chose représentée est une
convention: la lettre nommant un objet géométrique est un exemple de
symbole. L'\emph{indice} est un signe qui entretient un rapport physique ou
matériel avec la chose représentée: le mouvement de rotation des ailes
d'un moulin est un indice du vent. L'\emph{icône} est un signe qui a un
rapport de ressemblance avec la chose qu'il représente: le dessin d'un
arbre peut être considéré comme une icône de la notion d'arbre.

Dans la pratique de la géométrie depuis les Grecs, le diagramme tracé
explicitement sur un support n'est pas supposé être le véritable objet
de l'étude. La meilleure preuve en est que des imperfections dans le
tracé sont acceptées, parfois même conseillées. Bref, puisque ce
diagramme tient lieu de la configuration géométrique idéale sur
laquelle porte le raisonnement tout en étant non confondu avec elle,
il est un signe au sens de Peirce. Il n'est pas difficile d'identifier
la catégorie de ce signe puisque le symbole et l'indice sont d'emblée
disqualifiés. Le diagramme géométrique est donc une icône: sa vertu
est de représenter une configuration géométrique par un rapport de
ressemblance. D'après Peirce, cette ressemblance porte sur la forme de
la configuration, sachant que pour lui, \emph{la forme d'une chose est
  la structure de cette chose, c'est-à-dire l'ensemble des relations
  entre les parties constitutives de cette chose.} L'intérêt du
diagramme est qu'il rend perceptible cette forme abstraite par
l'observation directe, au même titre que quiconque peut percevoir une
configuration d'objets concrets située dans son champ visuel. Il est
notable que de nombreux secteurs des mathématiques utilisent
explicitement des diagrammes se présentant sous la forme de figures
géométriques (au sens d'une disposition concrète de points, de droites
et autres courbes, de lettres et de chiffres), alors même que les
objets étudiés ne relèvent pas de la géométrie. Pour ne donner qu'un
exemple emblématique, il suffit de citer la théorie des catégories
dans laquelle les diagrammes sont des figures constituées
principalement de flèches, et dont la géométrie représente des
propriétés abstraites sans rapport immédiat avec la géométrie au sens
usuel. La figure~\ref{fig:2} montre un diagramme de ce type illustrant
la notion de \emph{produit fibré} dans une catégorie.
\begin{figure}[h]
  \[
    \xymatrix@C+2em@R+2em{
      & N \ar@/_10pt/[ddl]_*{g''} \ar@/^10pt/[ddr]^*{f''} \ar@{.>}[d]|*+<5pt,5pt>{\alpha} & \\
      & P \ar[dl]_(0.4)*-<1ex,1ex>{g'} \ar[dr]^(0.4)*-<1ex,1ex>{f'} & \\
      A \ar[r]^*{f} & C & B \ar[l]_*{g} }
  \]
  \caption{Le diagramme de la propriété universelle du produit
    fibré~$(P,g',f')$ de deux morphismes~$f\colon A\to C$
    et~$g\colon B\to C$.}
  \label{fig:2}
\end{figure}
Toutes ces formes géométrisées de diagrammes relèvent clairement de
l'analyse précédente.  Poser qu'un diagramme en ce sens est une icône
donnant à voir la forme d'une configuration abstraite sur une figure
géométrique est une clarification intéressante, mais cela ne constitue
peut-être pas une grande surprise et encore moins une révolution dans
la philosophie des mathématiques. De fait, la position de Peirce est
bien plus audacieuse comme la suite va le montrer.

Tout d'abord, il introduit une extension de la notion de diagramme qui
coupe le lien qui liait cette notion avec celle de figure
géométrique. En toute généralité, un diagramme est une icône qui met
en scène de manière perceptible la forme de la chose représentée par
cette icône\footnote{En plus du diagramme,
  l'analyse percienne distingue deux autres types
  d'icônes: l'\emph{image}, qui désigne un signe
  perçu comme la représentation analogique de quelque chose --~l'image,
  c'est «la ressemblance à l'état pur», le «degré zéro de l'icône»~--,
  et la \emph{métaphore}, qui consiste à remplacer une relation
  sémiotique par une nouvelle relation sémiotique sur la base d'une
  similarité non préexistante introduite à cette occasion \citep{Verhaegen1994} . Quel que
  soit l'intérêt de ces notions par ailleurs fort discutées dans la
  littérature \citep{Verhaegen1994}, \citep{Esquenazi1997}, \citep{Lefebvre2013}, c'est le
  diagramme dans l'acception généralisée précédente qui constitue pour
  Peirce le processus sémiotique essentiel à l'{\oe}uvre dans les
  mathématiques. C'est cette thèse qui est retenue dans le présent
  article, et elle constitue l'un des appuis à la conception empirique
  de la démonstration donnée dans la partie~\ref{NonForm}.}.  Pour Peirce, il n'y a quasiment pas de limite quant à la
nature d'une entité pour qu'elle puisse être interprétée comme l'icône
d'une autre chose; la seule et impérieuse condition est que cette
entité ressemble à cette chose sous un certain rapport. Ce degré de
généralité se reporte pleinement sur la notion de diagramme comme la
suite va le montrer.

Comme exemple très simple de diagramme en ce sens, Peirce propose
l'arrangement suivant susceptible d'apparaitre dans une présentation
des bases de sa sémiotique:
\[
  \text{Signe:} \left\{
    \begin{aligned}
      &\text{symbole} \\
      &\text{indice} \\
      &\text{icône}
    \end{aligned}
  \right.
\]
En effet, il s'agit d'un dispositif graphique et textuel --~une
inscription scripto-visuelle~-- qui illustre de manière empiriquement
perceptible la propriété abstraite selon laquelle la classe des signes
est constituée de la réunion des classes des symboles, des indices et
des icônes.

Sa notion généralisée de diagramme permet à Peirce d'énoncer sa thèse
fondamentale sur la sémiologie des mathématiques: \emph{dans tous ses
  champs d'étude, c'est l'essence même de la pensée mathématique en
  acte de se traduire constamment par un travail de construction,
  d'observation et de transformation de diagrammes.} Cette thèse prend
le contre-pied de l'idée reçue selon laquelle le mathématicien accède
directement aux formes ou concepts abstraits; il est éventuellement
ajouté que cela résulte d'une faculté plus ou moins mystérieuse nommée
l'intuition. Il se trouve que le point de vue de Peirce englobe aussi
une interprétation éclairante de l'intuition mathématique parfaitement en accord
avec la précédente thèse, puisqu'il en donne la définition suivante:
\emph{l'intuition d'une forme abstraite est la saisie empirique de
  cette forme sur une icône qui la représente.}

Pour comprendre la pertinence de ce point de vue tout à fait radical,
on peut le tester sur le cas des formules algébriques qui peuplent les
traités d'algèbre. En voici un exemple extrait au hasard des
\emph{Éléments de mathématique} de Bourbaki \cite[page
A.I.93]{bourbaki70}.
\[
  \biggl( \sum_{\lambda \in\mathrm L}x_{\lambda}\biggr)^n = \sum_{|
    \beta |=n} \frac{n!}{\prod\limits_{\lambda \in\mathrm L}
    \beta_{\lambda}!}  \prod_{\lambda \in\mathrm
    L}x_{\lambda}^{\beta_{\lambda}}\text.\leqno(14)
\]
Pour Peirce, une formule de ce type, avec ses notations si
particulières et sa manière d'occuper l'espace sur le support de
l'écriture --~savoir-faire inventé graduellement au cours du
développement historique de ce domaine~-- est aussi un diagramme. En
effet, cette graphie donne à voir au mathématicien qui l'observe
nombre de propriétés de l'objet abstrait ainsi représenté, et prépare
les opérations qu'on pourra lui appliquer. Autre exemple, comment ne
pas être saisi par l'aspect diagrammatique de l'ensemble suivant de
formules?
\[
  \begin{aligned}
    (a+b)^1&=a+b \\
    (a+b)^2&=a^2+2ab+b^2 \\
    (a+b)^3&=a^3+3a^2b+3ab^2+b^3 \\
    \multicolumn{2}{c}{\dotfill}
  \end{aligned}
\]
C'est ainsi qu'il est possible d'interpréter une bonne partie de la
pratique de l'algèbre comme un travail d'observation et de
transformation de diagrammes. Mais bien entendu, celui qui se livre à
cette activité algébrique n'a pas à être conscient de l'interprétation
diagrammatique de son travail: il «fait de l'algèbre», un point c'est
tout!  Ce qui vient d'être dit sur les formules algébriques est
évidemment valable pour toute formule mathématique quel que soit le
domaine mathématique considéré.

Si on prend au sérieux le degré de généralité avec lequel quelque
chose est susceptible de jouer le rôle d'icône, qui pourrait être
n'importe quoi d'après Peirce, un fragment de langage lui-même peut
être porteur d'une dimension diagrammatique. Précisément, une séquence
purement textuelle dans un développement mathématique constitue un
diagramme, sous réserve qu'elle se révèle susceptible de mettre en
scène une forme abstraite de manière perceptible. Dans cet ordre
d'idée, la structure en «formules» du texte mathématique grec mise en
lumière par Netz est interprétable comme l'explicitation de la
dimension diagrammatique de ce texte, dimension indispensable à la
fécondité de l'activité mathématique.  On peut ajouter qu'en général,
séquences textuelles, formules mathématiques et éventuellement figures
géométriques sont inséparables et concourent à la constitution de
diagrammes composés qui sont des mixtes de ces divers
ingrédients. D'ailleurs, le diagramme géométrique introduit par les
Grecs est lui même un diagramme composé puisqu'il est constitué de
lettres et de dessins.  Enfin, une entité mentale peut être un
diagramme, si cette entité est ressentie comme satisfaisant la
définition de Peirce. Finalement, les diagrammes au sens général
apparaissent comme le substrat de l'activité mathématique. Cette
activité se réalise par diverses opérations portant sur des
diagrammes, comme les constructions et transformations de
diagrammes. D'après Peirce, les transformations de diagrammes
apportent la preuve visible que les mathématiques s'enrichissent par
la production de nouvelles connaissances puisque en effet, ces
opérations ne sont nullement prédéterminées.

\subsection{Le point sur la question d'un mode d'existence en
  mathématiques après Netz et Peirce}

Les travaux de Netz et de Peirce présentés dans les parties
précédentes portent tous les deux sur les mathématiques, mais
diffèrent notablement, tant par le champ d'étude --~mathématiques
grecques et mathématiques en général~-- que par le type d'approche
--~méthodes empiriques d'analyse textuelle des sciences sociales et
analyse philosophique et sémiotique~--, au point qu'il peut sembler
aventureux de vouloir en tirer une leçon commune. Et pourtant, ces
mêmes différences font que ces deux approches des mathématiques se
complètent relativement bien en ce qui concerne la place du diagramme
dans les mathématiques. En effet, le diagramme géométrique avec
lettres dont Netz a analysé l'importance dans les mathématiques
grecques est généralisé par Peirce en un dispositif iconique jouant un
rôle essentiel dans l'ensemble des mathématiques et non plus dans la
seule géométrie grecque. D'autre part, Netz a mis en évidence des
facteurs linguistiques spécifiques chez les Grecs --~utilisation d'un
langage au lexique réduit et structuration de ce langage par l'emploi
répétitif de formules~--, facteurs dont on peut admettre la
persistance au delà de l'Antiquité jusqu'à nos jours et qui
constituent un cadre dans lequel la notion de diagramme généralisé au
sens de Peirce se révèle pertinente pour l'activité mathématique en
général. Il s'agit maintenant d'analyser dans quelle mesure cette
leçon commune permet de progresser relativement à la question d'un
mode d'existence pour les mathématiques au sein de l'\textsc{eme}.

Préalablement, il est intéressant de noter que l'un des bénéfices
majeurs de l'introduction dans l'\textsc{eme} du mode [\textsc{ref}]
de la référence objective est que ce mode n'a pas vocation à être
l'unique vecteur de toute l'ontologie. Au contraire, cette conception
du travail scientifique laisse de la place à d'éventuels autres modes
d'existence capables de prendre en charge de nouveaux êtres, de les
faire exister, de les instaurer. C'est effectivement ce qui est fait
dans \citep{Latour2012,EME} par l'introduction d'une série de modes
d'existence pour, par exemple, les êtres de la fiction, les êtres de
la technique, les êtres du droit, de la politique, de la parole
religieuse, etc. Chacun de ces nouveaux modes d'existence se manifeste
par le truchement de trajectoires (continuité) constituées de sauts
entre des formes différentes (discontinuité) se succédant selon un
principe de «véridiction» \cite[«Livre», colonne «Vocabulaire», entrée
«Véridiction»]{EME} spécifique à chaque mode.  La richesse de ce
nouveau point de vue sur l'ontologie, sa capacité à rendre compte de
pratiques aussi éloignées que par exemple celles relevant d'une part
de la recherche scientifique et d'autre part de la parole religieuse,
laisse penser qu'il devrait pouvoir permettre une approche originale
et pertinente du cas des mathématiques.

Les analyses de Netz et de Peirce constituent un premier pas dans
cette direction.  En effet, la description qu'ils donnent des
mathématiques est manifestement en forte corrélation avec le fond
philosophique qui initie l'analyse anthropologique dans la démarche
scientifique développée par l'\textsc{eme}. C'est vrai en ce qui
concerne la critique vis-à-vis de l'approche épistémologique qui
privilégie les révolutions scientifiques au détriment d'une étude
concrète du fonctionnement réel d'une science normale. C'est aussi
vrai pour l'importance donnée aux véritables technologies
intellectuelles que sont d'un côté le travail sur les diagrammes
couplé avec un langage enrégimenté, de l'autre la fabrication et
l'étude des inscriptions scripto-visuelles. Dans les deux cas, il
s'agit de «rematérialiser» la connaissance mathématique, c'est-à-dire
de montrer en quoi elle est dépendante de procédés concrets de
construction et de manipulation d'inscriptions, dont l'importance est
généralement ignorée par l'approche purement épistémique du savoir
mathématique.

Néanmoins, ce que l'on ne trouve pas dans les analyses de Netz et de
Peirce, ce qui d'ailleurs n'est pas leur propos, c'est la mise en
évidence des cheminements particuliers --~l'analogue de la notion de
chaine de référence dans les sciences expérimentales~-- par lesquels
se manifesterait l'existence des entités mathématiques. Toutefois,
Netz analyse très précisément la forme des démonstrations des Grecs.

Après Netz et Peirce, que peut-on dire alors sur la manière dont se
manifestent les êtres mathématiques? Ce sont les «existants»
\citep[«Livre», colonne «Vocabulaire», entrée «Êtres, existants»]{EME}
qui résultent des agencements particuliers produits par les
mathématiciens.  Ces agencements sont faits de textes, de dessins et
de diagrammes constituant des ensembles éminemment
structurés. Essentiellement, on y voit, et cela depuis les anciens
Grecs, se répéter la forme suivante: l'introduction d'un cadre de
travail par un jeu de définitions, de diagrammes, de notations et de
commentaires\footnote{Dans les \emph{Éléments} d'Euclide, cette
  introduction n'est accompagnée d'aucun commentaire. D'autres traités
  antiques contiennent une préface, mais celle-ci n'a en général pas
  de finalité méthodologique.}, puis l'annonce du résultat à démontrer
suivie de sa validation par une démonstration proprement dite. La
leçon de Netz et de Peirce est que ces formes ne deviennent
productives et partageables par une communauté que sous certaines
conditions techniques et cognitives qui caractérisent ces
mathématiques, précisées dans les parties~\ref{emerge} et~\ref{Peirce}
ci-dessus.  Marquées par l'usage des démonstrations comme moyen de
validation ultime des résultats, ce sont ces mathématiques qui,
héritières des Grecs, se sont trouvées au cœur de la révolution
scientifique et sont maintenant pratiquées dans les lieux de recherche
partout dans le monde. Dans la présente étude, ce sont elles qui sont
désignées par l'expression «les mathématiques\footnote{L'expression
  est trompeuse car elle peut donner à penser que «les» mathématiques
  en question sont parfaitement unifiées et qu'elles épuisent la
  totalité des formes de l'activité mathématique. Cela est tout à fait
  discutable, ne serait-ce qu'au niveau historique. Un exemple
  frappant d'une autre pratique mathématique est donné par l'étude
  \citep{Chemla2016} de textes datant de la Chine impériale du~3\ieme\
  au 7\ieme~siècle. On y découvre une culture mathématique riche mais
  très différente de celle qui est au centre de l'attention du présent
  travail. Malgré cette altérité, certaines des caractéristiques des
  mathématiques grecques relevées par Netz semblent présentes, en
  particulier un usage, en partie caché, d'inscriptions
  scripto-visuelles.}» et qui constituent son objet central.

Ceci étant admis, la place centrale des démonstrations dans les
mathématiques\footnote{À juste titre, les mathématiciens n'aiment
  guère que l'on réduise leur activité à l'écriture de
  démonstrations. De fait, une démonstration en bonne et due forme
  n'est susceptible d'être rédigée qu'après une longue phase
  préparatoire relativement informelle. Il n'empêche que c'est
  l'élaboration finale d'une démonstration qui est le but de ce
  travail préalable et qui en marque le succès. De plus, la phase
  préparatoire informelle est difficile à observer puisqu'elle
  comporte souvent des phases de réflexion intense sans manifestation
  extérieure significative. On peut écouter la belle description qu'en
  donne André Lichnerowicz dans l'interview \citep{vernier87} dont la
  pointe est le propos rapporté que «le plus difficile pour une femme
  de mathématicien, c'est de distinguer un mathématicien qui travaille
  et un mathématicien qui dort».} suggère que les êtres mathématiques
manifestent leur existence par des parcours en relation avec les
démonstrations dans lesquelles ils sont impliqués. Est-il possible
d'interpréter sans distorsion manifeste ces parcours démonstratifs
comme des chaines de référence au sens de l'\textsc{eme}? Ou bien,
peut-on concevoir les démonstrations comme les trajectoires d'un
nouveau mode d'existence, qui lui serait spécifique aux entités
mathématiques? Pour répondre à ces questions, un long détour pour
élaborer une conception non formaliste des mathématiques, conception
venant elle-même en critique de la vision formaliste de cette
discipline, s'avère utile.

\section{La démonstration et l'idée formaliste de sa
  perfection} \label{StrucPreuve}
 
La logique contemporaine met à disposition plusieurs types de
description de la notion de démonstration captant plus ou moins bien
la forme de la pratique argumentaire des mathématiciens. Le calcul
nommé «déduction naturelle» est particulièrement
intéressant pour ce propos et la description qui suit s'en inspire un
peu, en évitant toutefois les aspects techniques moins pertinents dans
ce travail philosophique. Introduite par le mathématicien et logicien
allemand Gerhard Gentzen\footnote{Gentzen (1909-1945) présente la
  curiosité d'avoir été à la fois un très grand logicien et un nazi
  semble-t-il convaincu. Il serait malvenu d'en déduire que son œuvre
  en logique --~la déduction naturelle, le calcul des séquents, la
  cohérence de l'arithmétique de Peano~-- présente une quelconque
  connivence avec le national-socialisme. De fait, son apport
  scientifique s'est révélé d'une extrême importance pour les
  développements modernes de la théorie de la démonstration.} \citep{Gentzen1935}, la
déduction naturelle a pour objet premier de proposer un formalisme
reproduisant le mode naturel du raisonnement mathématique par
opposition au style de la méthode axiomatique\footnote{La méthode
  axiomatique n'est pas représentative de la pratique courante des
  mathématiciens. La meilleure preuve en est que la grande majorité
  d'entre eux n'ont que de très vagues idées sur le système
  axiomatique qui est censé être le fondement de leur travail. Cela ne
  les empêche aucunement de faire des mathématiques!} commune aux
systèmes de Frege, de Russell et Whitehead et de Hilbert. Selon
Gentzen, le raisonnement mathématique naturel repose non pas sur la
référence à une poignée d'axiomes trop souvent abscons mais sur
l'application de règles d'inférences claires, efficaces et en aussi
grand nombre que nécessaire dans chaque secteur de cette discipline.

Selon la terminologie utilisée dans les travaux de Per
Martin-Löf \citep{MartinLof87}, un résultat mathématique dont on sait
qu'on a le droit de l'affirmer est nommé un \emph{jugement}, ce terme
étant pris exclusivement dans le sens de \emph{jugement mathématique}
dans toute la suite du présent texte. Un jugement est donc un élément
établi de connaissance mathématique: il se présente sous la forme d'un
énoncé mathématique qui a été préalablement reconnu comme valide.

La forme la plus simple d'un jugement est «$A$~est vrai\footnote{Le
  jugement «$A$ est vrai» est synonyme de «$A$ est un théorème» dans
  certains écrits de logique. Néanmoins, la tradition mathématique
  veut que l'on réserve le terme de théorème aux résultats qui sont à
  la fois importants et difficiles à démontrer.  Il y a d'autres
  formes de jugement, en particulier le jugement hypothétique «$A$ est
  vrai sous l'hypothèse que $B$ est vrai».}», où $A$~désigne une
proposition. Par exemple, l'énoncé que «la somme des angles d'un
triangle est égale à deux angles droits» est une proposition, et
l'affirmation que cette proposition est vraie est un jugement qui est
par exemple établi par le raisonnement que l'on trouve dans les
\emph{Éléments} d'Euclide et que nous présentons dans la
partie~\ref{PreuveExpe}.

Maintenant, il reste à savoir de quel droit on peut affirmer un jugement.
Certains sont donnés préalablement à la phase déductive; ils
traduisent une sorte d'évidence basique et immédiate relative aux
objets considérés et au cadre de travail.  Les autres jugements sont
tous ceux que l'on peut obtenir par déduction, c'est-à-dire par
application de \emph{règles d'inférence}: chacune de ces règles décrit
comment zéro, un ou plusieurs jugements permettent d'établir un autre
jugement.  Ces règles d'inférence reprennent sous une forme épurée les
bases de la logique du raisonnement humain. La déduction naturelle
représente une démonstration mathématique comme un enchainement de
jugements\footnote{\label{StructDem}Plus généralement, la structure
  d'une démonstration peut être un peu plus complexe qu'une simple
  suite organisée comme un parcours linéaire: c'est plutôt une famille
  de jugements reliés par des règles d'inférence qui confluent vers le
  jugement final. La structure combinatoire d'une telle configuration
  est ce que les logiciens nomment un \emph{arbre de preuve}. Cette
  même structure serait vraisemblablement appropriée pour les chaines
  de référence. Dans le présent travail de nature essentiellement
  philosophique, il semble inutile de préciser davantage ce point
  technique.}, chaque jugement étant établi à partir du précédent (ou
des précédents) par une règle d'inférence, le dernier jugement établi étant l'affirmation du théorème.

\emph{A priori}, il n'y a pas grand-chose à reprocher à la description
de la notion de démonstration développée ci-dessus sinon qu'elle est
excessivement brève et formelle. Toutefois, ce mode de description
donne facilement prise à une interprétation formaliste, cette dernière
se révélant être une version moderne de l'antique captation
philosophique des mathématiques grecques telle qu'elle avait été
opérée par les philosophes platoniciens.

En effet, on a vu qu'une démonstration apparait comme une structure
constituée de jugements reliés par des inférences logiques; or les
inférences logiques sont parfaitement caractérisées et cataloguées par
le travail des logiciens depuis le siècle dernier: cette structure
peut donc être simulée par un calcul en réduisant le texte de la
démonstration à une série de symboles assemblés selon des règles
syntaxiques précises, sans utilisation de la langue naturelle. Il en
découle que la démonstration devient un dispositif formel objectif
dont la validité tient à son objectivité, c'est-à-dire à son caractère
indépendant de l'approche subjective de l'esprit humain. Autrement
dit, la validité d'une démonstration réside ultimement dans une
propriété objective d'un objet combinatoire extérieur à notre
conscience. En conséquence, une fois qu'elle est convenablement
formalisée, une démonstration mathématique constitue une forme
d'argumentation imparable, incontestable, menant inéluctablement à sa
conclusion. Ainsi, la vérité dévoilée par une démonstration acquiert un
caractère universel, immuable et permanent, bien exprimé par la
conception tarskienne de vérité.

Afin d'illustrer la réalité de ce point de vue, voici une longue
citation issue de l'introduction du livre~I des \emph{Éléments de
  mathématique} de Bourbaki \citep{bourbaki54}, intitulé \emph{Théorie
  des ensembles}, qui a pour but explicite de montrer les changements
--~par rapport au vénérable héritage légué par les Grecs~-- que les
développements contemporains apportent à la notion de démonstration.
\begin{quotation}\label{citation-bourbaki}
  [\dots] l'analyse du mécanisme des démonstrations dans des textes
  mathématiques bien choisis a permis d'en dégager la structure, du
  double point de vue du vocabulaire et de la syntaxe. On arrive ainsi
  à la conclusion qu'un texte mathématique suffisamment explicite
  pourrait être exprimé dans une langue conventionnelle ne comportant
  qu'un petit nombre de «mots» invariables assemblés suivant une
  syntaxe qui consisterait en un petit nombre de règles inviolables:
  un tel texte est dit \emph{formalisé}. La description d'une partie
  d'échecs au moyen de la notation usuelle, une table de logarithmes,
  sont des textes formalisés; les formules du calcul algébrique
  ordinaire en seraient aussi, si l'on avait complètement codifié les
  règles gouvernant l'emploi des parenthèses et qu'on s'y conformât
  strictement, alors qu'en fait certaines de ces règles ne
  s'apprennent guère qu'à l'usage, et que l'usage autorise à y faire
  certaines dérogations.

  La vérification d'un texte formalisé ne demande qu'une attention en
  quelque sorte mécanique, les seules causes d'erreur possibles étant
  dues à la longueur ou à la complication du texte; c'est pourquoi un
  mathématicien fait le plus souvent confiance à un confrère qui lui
  transmet le résultat d'un calcul algébrique, pour peu qu'il sache
  que ce calcul n'est pas trop long et a été fait avec soin. Par
  contre, dans un texte non formalisé, on est exposé aux fautes de
  raisonnement que risquent d'entrainer, par exemple, l'usage abusif
  de l'intuition, ou le raisonnement par analogie. En fait, le
  mathématicien qui désire s'assurer de la parfaite correction, ou,
  comme on dit, de la «rigueur» d'une démonstration ou d'une théorie,
  ne recourt guère à l'une des formalisations complètes dont on
  dispose aujourd'hui, ni même le plus souvent aux formalisations
  partielles et incomplètes fournies par le calcul algébrique et
  d'autres similaires; il se contente en général d'amener l'exposé à
  un point où son expérience et son flair de mathématicien lui
  enseignent que la traduction en langage formalisé ne serait plus
  qu'un exercice de patience (sans doute fort pénible). Si, comme il
  arrive mainte et mainte fois, des doutes viennent à s'élever, c'est
  en définitive sur la possibilité d'aboutir sans ambigüité à une
  telle formalisation qu'ils portent, soit qu'un même mot soit employé
  en des sens variables suivant le contexte, soit que les règles de la
  syntaxe aient été violées par l'emploi inconscient de modes de
  raisonnement non spécifiquement autorisés par elles, soit encore
  qu'une erreur matérielle ait été commise.  Ce dernier cas mis à
  part, le redressement se fait invariablement, tôt ou tard, par la
  rédaction de textes se rapprochant de plus en plus d'un texte
  formalisé, jusqu'à ce que, de l'avis général des mathématiciens, il
  soit devenu superflu de pousser ce travail plus loin; autrement dit,
  c'est par une comparaison, plus ou moins explicite, avec les règles
  d'un langage formalisé, que se fait l'essai de la correction d'un
  texte mathématique.
\end{quotation}

C'est ainsi que la mise au clair de ce qu'est une démonstration au
moyen d'un point de vue avant tout formaliste peut conduire à
conforter l'idée de la perfection absolue de l'argumentation
mathématique. Non seulement cette conception mérite d'être
questionnée, mais, par un glissement insidieux, elle fait de la
démonstration mathématique le modèle par excellence de la méthode
rationnelle permettant, selon «les modernes\footnote{Conformément à
  l'utilisation de ce terme dans l'\textsc{eme} \citep[«Livre», colonne
  «Vocabulaire», entrée «Modernes, modernisation»]{EME}, «les
  modernes» désignent ces personnes qui adhèrent pleinement au grand
  récit de l'émancipation selon lequel, depuis la révolution
  scientifique, s'établirait progressivement mais irréversiblement la
  distinction absolue entre Illusion et Raison, entre un passé
  archaïque et un futur émancipé, entre Croyance et Savoir Vrai. Voir
  l'\textsc{eme} \citep{Latour2012,EME} pour de plus amples
  développements sur ce thème et son rapport avec la crise
  écologique.}», l'accès direct et complet à la Vérité quel que soit
le domaine ou le problème considéré.

\section{Une conception empirique du jugement
  mathématique} \label{NonForm}
 
\subsection{L'expérience d'un jugement} \label{PreuveExpe}

L'interprétation formaliste précédente peut induire le biais qu'un
jugement serait obtenu en tant que tel grâce au dispositif formel que
constitue une démonstration formalisée, indépendamment des acteurs
mathématiciens concernés, de leur appréhension de ladite
démonstration, de leurs tentatives, de leurs discussions, de leurs
calculs. Pour échapper à ce travers, il suffit de remarquer qu'une
configuration de signes sur une feuille de papier ou un tableau ne
devient éventuellement une démonstration --~formelle ou pas~-- qu'à la
condition que des acteurs humains puissent en faire
l'expérience. Cette remarque nous met sur la voie d'une autre approche
philosophique, cette fois-ci dans l'esprit de l'\emph{empirisme
  radical} de William James \citep{James2007}. Cette doctrine est
centrée sur la notion d'expérience, ce dernier terme se référant à
tout ce dont un sujet individuel fait l'expérience dans le flux
continu de sa propre vie. Autrement dit, l'expérience est le champ
instantané du présent, à la fois continu et changeant, qui constitue
le vécu de chaque personne. L'intérêt de
l'empirisme jamesien pour notre questionnement se manifeste dans le principe suivant qui en
offre un résumé percutant:

\begin{quote}
  \emph{N'admettre que ce dont on puisse faire l'expérience, et rendre
    justice à tout ce qui peut être objet d'expérience}\footnote{Cette
    formulation provient de la préface de \citep{James2007}.}.
\end{quote}

De ce point de vue, un jugement, plus exactement l'affirmation d'un
jugement ou l'adhésion à un jugement, est une expérience d'un certain
type qu'il faut maintenant essayer de décrire.
 
La question est de savoir ce qui caractérise les jugements parmi les
expériences porteuses de connaissance. La réponse qui va être donnée
est inspirée par la lecture du travail \citep{MartinLof87} du logicien,
mathématicien et philosophe Per Martin-Löf, connu entre autres comme
le concepteur de l'étonnante et riche \emph{théorie constructive des
  types}. Néanmoins, les développements qui vont suivre ne sont pas
dans la droite ligne de l'analyse de Martin-Löf. De fait, cet auteur
ne se réfère pas directement à l'empirisme ou au pragmatisme de
James. Ses sources sont plutôt l'œuvre de Kant et surtout la
phénoménologie à travers les analyses de Husserl sur les fondements de
la logique\footnote{\label{IntuiML}De plus, dans ses travaux,
  Martin-Löf se réfère explicitement à l'intuitionnisme dont il est un
  représentant ancré dans la période contemporaine \emph{post}
  1960. Contrairement au maitre Brouwer initiateur de ce courant, il
  ne jette pas l'opprobre sur le langage et les langues
  formelles. Néanmoins, il est profondément fidèle à l'exigence de
  sens qui constitue le trait principal de l'intuitionnisme. Cela se
  retrouve dans son interprétation de la démonstration d'un jugement
  qu'il donne dans \citep{MartinLof87}, qui a inspiré la définition
  donnée dans le présent travail. Voir la présentation synthétique
  \citep{Bourdeau2004} pour plus de précisions sur l'intuitionnisme et
  l'apport de Martin-Löf.}. Néanmoins, l'intuitionnisme auquel se
réfère Martin-Löf n'est pas totalement étranger à l'empirisme, comme
relevé dans \citep{Bourdeau2004}.
 
La propriété fondamentale qui fait qu'un énoncé mathématique est un
jugement est qu'\emph{on sait qu'on a le droit de l'affirmer}: il
possède une démonstration. Pour l'instant, rien de bien nouveau. La
divergence avec le point de vue formaliste se marque si l'on insiste
sur la nature de cette connaissance simplement en mettant en avant la
notion d'expérience. C'est ce qui est fait avec la définition suivante
qui va jouer un rôle essentiel dans les développements qui suivent.
 
\begin{quote} \label{DefPreuve} \emph{Une \emph{expérience de démonstration}
    d'un jugement est une expérience par laquelle un mathématicien
    saisit l'évidence de ce jugement}.
\end{quote}

On peut dire qu'une expérience de démonstration est une expérience par
laquelle ce mathématicien littéralement \emph{voit} ce jugement;
alors, pour le sujet de cette expérience, le jugement \emph{est},
pleinement, sans restriction au point qu'il lui est impensable d'en
douter, de s'y soustraire. Chaque terme de la définition a son
importance: «expérience» bien entendu pour qualifier un segment du
flux de la vie propre à un individu (ce qui n'empêche pas que la
conclusion de l'expérience puisse être partagée: ce point sera discuté
dans la partie~\ref{FormeDem}), «mathématicien» pour indiquer que la
reconnaissance pleine et entière de cette expérience nécessite une
compétence généralement acquise par l'enseignement et la pratique au
sein d'une communauté scientifique spécifique, enfin «l'évidence» pour
marquer qu'une expérience de démonstration rend visible et manifeste
le jugement, qu'elle expose au regard la totalité d'un arrangement qui
le valide sans conteste et l'établit comme fait mathématique. Dans
l'utilisation qui vient d'en être faite, les expressions «voir un
jugement», «rendre visible» ou «exposer au regard» ne sont pas de
simples métaphores destinées à désigner la saisie directe par l'esprit
de propriétés abstraites. Comme l'a montré Peirce, la dimension
diagrammatique des mathématiques permet au mathématicien d'amener dans
sa perception sensible (en général la vision) la forme de la
configuration abstraite dont il traite, et il peut ainsi l'observer et
la manipuler avec la force brute que donne cette perception.

\begin{figure}[!h]
  \centering \includegraphics[width=8cm]{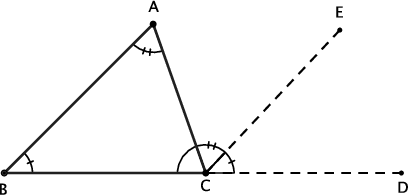}
  \caption{Un diagramme avec lettres analogue à celui qui accompagne
    la Proposition~32 du Livre~I des \emph{\'Eléments}
    d'Euclide.}\label{fig:3}
\end{figure}

Afin d'illustrer ce qu'est cette expérience de la vision d'un
jugement, considérons l'exemple de la Proposition~32 du Livre~I des
\emph{\'Eléments} d'Euclide qui énonce le jugement selon lequel
\emph{la somme des angles d'un triangle est égale à deux droits}. Cet
exemple est intéressant du fait qu'il est relativement élémentaire
tout en étant caractéristique de la démarche inaugurée par les
Grecs. La démonstration est basée sur le diagramme présenté dans la
figure~\ref{fig:3}. Sur ce diagramme se distinguent en traits continus
la donnée initiale d'un triangle~ABC et en traits discontinus une
construction auxiliaire telle que B,~C et~D soient alignés et se
succèdent dans cet ordre sur cette droite et que CE~soit parallèle
à~BA, le point~E étant situé du même côté de la droite prolongeant~BC
que le point~A.
Imaginons l'expérience vécue par un mathématicien qui scrute ce
diagramme. En se remémorant le schéma général des angles définis par
deux droites parallèles coupées par une sécante et les relations
d'égalité afférentes, cet observateur note les égalités d'angles
suivantes: $\widehat{\mathrm{ABC}}=\widehat{\mathrm{ECD}}$ et
$\widehat{\mathrm{BAC}}=\widehat{\mathrm{ACE}}$. Alors ce scrutateur
géomètre ne peut qu'être saisi par l'évidence du jugement considéré
qui se donne à voir complètement et clairement sur le diagramme.

Cette expérience qu'est une expérience de démonstration d'un jugement
est en partie un acte mental, mais elle ne se réduit pas au moment
ultime pendant lequel le mathématicien est saisi par cette
évidence. Ce moment est préparé par la démarche souvent laborieuse
d'analyse critique et d'appropriation à laquelle se livre le
mathématicien en composant ou en décryptant une inscription faite de
textes, de figures, de schémas, de formules, de calculs, etc. Ce n'est
qu'à la fin de cette phase d'exploration et de confrontation que cette
expérience peut éventuellement surgir pleine et entière.  C'est la
totalité du processus menant à l'évidence finale qui est l'expérience
de démonstration. Autrement dit, le succès de la phase finale de
l'expérience de démonstration se manifeste dans la conscience du
mathématicien mais il dépend grandement de l'inscription
scripto-visuelle sur laquelle s'appuie ce processus et plus
généralement des conditions cognitives caractéristiques des
mathématiques et de la dimension diagrammatique de cette discipline
mises en lumière dans les travaux de Netz et Peirce présentés dans la
partie~\ref{NetzPeirce}.

Contrairement à l'approche formaliste de la démonstration évoquée dans
la partie~\ref{StrucPreuve}, une expérience de démonstration d'un
jugement ne se fonde pas sur l'utilisation explicite de techniques et
méthodes logiques. Ce qui justifie une expérience de démonstration
d'un jugement par un mathématicien, c'est la mise en place par ce
dernier d'un dispositif dont l'expérience est vécue comme la
manifestation incontestable du bien-fondé --~l'évidence~-- de ce
jugement. Dans l'ordre conceptuel, cette manifestation de l'évidence
mathématique précède l'analyse logique du raisonnement, qui elle
procède d'une mise en forme postérieure à l'expérience de
démonstration. En général, cette évidence n'est pas immédiatement
accessible. Le travail de construction d'une expérience de
démonstration auquel se livre le mathématicien consiste à fabriquer un
cheminement souvent indirect, permettant de transporter étape par
étape l'évidence jusqu'à obtenir l'évidence finale du jugement. Ce
transport de l'évidence produit à chaque étape une expérience de
démonstration d'un jugement intermédiaire. Ainsi, il apparait qu'une
expérience de démonstration d'un jugement se scinde elle-même en une
succession d'expériences de jugements intermédiaires, c'est-à-dire en
une succession d'expériences de l'évidence de chacun de ces jugements.  Nous proposons de lire dans ce sens le passage du \emph{Continu} suivant de Hermann Weyl.
\begin{quotation}
  Le célèbre ouvrage de Dedekind, \emph{Nature et fonction du nombre},
  1888, préf.\ 1\iere\ éd., s'ouvre sur cette phrase: «Ce qui est
  démontrable ne doit pas, en science, être cru sans
  démonstration». Cette assertion est probablement typique de la
  pensée de la plupart des mathématiciens, elle n'en prend pas moins
  les choses à l'envers. Comme si un rapport de fondation aussi médiat
  que celui que nous nommons démonstration, était capable de susciter
  une «croyance» sans que nous nous assurions de l'exactitude de
  chaque pas isolé au moyen d'une intuition immédiate! C'est cette
  intuition (non pas la démonstration), qui est partout la source
  ultime de la connaissance, elle consiste à avoir «l'expérience
  directe du vrai\footnote{Weyl cite ici Edmund Husserl: «Une évidence
    n'est, au contraire, rien d'autre que ``l'expérience vécue'' de la
    vérité» \citep[§~51]{husserl59}.}». \citep[p.~59.]{weyl94}
\end{quotation}
Ainsi, par l'expérience de démonstration globale, l'évidence se
propage de jugement en jugement pour culminer par l'évidence du
jugement final.

Mais d'où vient la puissance de la conviction attachée à une
expérience de démonstration d'un jugement, cette certitude qui a
impressionné les philosophes au point d'y voir le paradigme d'un accès
direct à la Vérité?  Pourquoi la «certitude nécessaire», c'est-à-dire
le type de confiance que l'on accorde à un résultat qui résulte d'une
expérience de démonstration, est-elle considérée comme l'une des
formes de certitude la plus incontestable qui soit?  Il y aurait lieu
d'approfondir ces questions et de relativiser certaines de ces
formulations, par exemple en discutant les arguments développés par
Wittgenstein dans \citep{Wittgenstein2006} et par Dewey dans
\citep{Dewey2014}, ce qui n'est pas fait sinon de manière allusive dans
le cadre de la présente contribution.  Pour l'instant, il suffit de
constater que c'est la notion même d'expérience de démonstration telle
qu'elle vient d'être définie qui explique la puissance de la certitude
nécessaire. Celui qui éprouve une expérience de démonstration d'un
jugement, c'est-à-dire celui qui est saisi par l'évidence du résultat
énoncé, est absolument certain de la validité du jugement considéré. À
l'opposé, celui qui admet qu'un résultat est démontré mais dont la
conviction est mitigée d'une dose de doute, celui-ci n'a simplement
pas éprouvé une expérience de démonstration de ce résultat. Comme le
montrera le point~2 de la partie~\ref{DeuxRem}, ce n'est pas une
question de force qui distingue la conviction attachée à une
expérience de démonstration d'un jugement, c'est une question de
qualité, celle qui résulte du cheminement de l'évidence qui aboutit en
pleine lumière, sans zone d'ombre et sans délégation à un
oracle\footnote{\label{NoteOracle}Un oracle contemporain pourrait
  être: l'avis d'un expert assurant d'autorité que telle démonstration
  est valide, ou bien encore le recours à un calcul sur machine
  approprié pour gérer tout ou partie d'une démonstration (voir le
  point~2 de la partie~\ref{DeuxRem}). En toute rigueur, l'argument
  selon lequel tel jugement est acquis du seul fait qu'il a fait
  l'objet d'une publication est aussi une forme de recours à un
  oracle\dots}, au jugement considéré. Rappelons encore une fois que
la nature même d'une expérience de démonstration et la qualité de la
certitude nécessaire qui en découle, sont profondément dépendantes des
conditions et dispositifs propres aux mathématiques précisés par Netz
et Peirce\footnote{\label{note:voevodsky}D'ailleurs, c'est lorsque ces
  conditions sont mises à mal par l'étendue et la complexité du champ
  mathématique considéré que la certitude nécessaire peut commencer à
  être mise en doute; c'est par exemple le sens de certaines
  déclarations du mathématicien Vladimir Voevodsky à propos de son
  absence de confiance en de nombreux résultats célèbres et publiés de
  son domaine de recherche, y compris les siens.}.

\subsection{Inscription de démonstration et partage d'un
  jugement}\label{FormeDem}

Maintenant qu'a été introduit le concept d'expérience de
démonstration, il est notable que cette notion ne recouvre pas la
totalité des significations contenues dans l'idée usuelle de
démonstration. En effet, dans un texte de mathématiques (livre,
article ou lettre), la partie du document qui est censée justifier un
théorème est elle-même qualifiée de démonstration (ou de preuve) du
théorème. Or il est clair que cet objet concret qu'est la
démonstration en ce sens est quelque chose de différent de
l'expérience de démonstration du théorème en question. C'est d'autant
plus clair pour un mathématicien qui sait que, généralement, la simple
lecture d'une démonstration d'un théorème est loin de suffire pour en
éprouver l'évidence. Sa pratique des mathématiques lui a montré
qu'avant d'atteindre cette évidence, il est souvent nécessaire de se
livrer à un travail d'exploration de cette démonstration, et que cette
tâche est délicate, et son issue incertaine. La définition qui suit
permet de préciser utilement cette distinction.

\begin{quote}
  \emph{Une \emph{inscription de démonstration} d'un jugement donné est une
    inscription scripto-visuelle (couchée sur le papier, tracée sur un
    tableau, enregistrée dans un fichier informatique\dots) présentant
    une configuration d'arguments textuels et graphiques susceptible
    en principe de permettre ou de favoriser l'expérience de
    démonstration du jugement.}
\end{quote}

La notion d'inscription de démonstration d'un jugement est donc
intrinsèquement différente de celle d'expérience de démonstration
puisque la première est une forme faite de textes, de tableaux, de
schémas, de diagrammes, et que la seconde est une expérience se
concluant chez celui qui l'éprouve par l'acte mental de l'évidence du
jugement.

Néanmoins, ces deux notions sont liées. Ou bien l'inscription de la
démonstration est élaborée en même temps que l'auteur cherche à
éprouver la validité d'un énoncé, validité qui est finalement attestée
par l'expérience de la démonstration. Dans ce cas, la recherche de
cette expérience oriente la fabrication d'une inscription, et cette
dernière, au fur et à mesure de son élaboration, offre un outil
scripto-visuel stable facilitant l'avènement de l'expérience de la
démonstration. Ou bien on dispose \emph{a priori} d'une démonstration
déjà rédigée, et un travail d'investigation basé sur cette inscription
est susceptible d'être constitutif d'une expérience de démonstration
du jugement considéré. Ainsi, \emph{une inscription de démonstration
  d'un jugement joue le rôle d'un guide pour le cheminement de
  l'évidence jusqu'à l'expérience de démonstration de ce
  jugement}\footnote{C'est un peu le même rapport qui existe entre le
  topo-guide décrivant une randonnée (bel exemple d'inscription
  scripto-visuelle) et l'accomplissement sur le terrain de cette même
  randonnée (expérience éventuellement délicate et
  incertaine).}. Finalement, c'est l'idée habituelle de démonstration
qui est scindée en les deux notions différentes mais complémentaires
que sont l'inscription de démonstration et l'expérience de
démonstration\footnote{Cette manière de distinguer deux notions
  différentes dans celle de démonstration est d'usage purement
  philosophique. Elle ne vise aucunement à s'imposer dans la pratique
  des mathématiciens.}.

L'intérêt principal d'une inscription de démonstration d'un jugement
est de conserver et de communiquer des éléments permettant de valider
ce jugement. L'acte de communiquer la démonstration d'un jugement est
essentiel dans la vie scientifique de la communauté des
mathématiciens. Cet acte présuppose la croyance que la démontrabilité
d'un jugement est une notion partageable.

Pourtant, il va de soi qu'une expérience de démonstration d'un
jugement est essentiellement privée du fait que l'expérience éprouvée
par un individu lui appartient en tant que moment de sa vie propre et
que toute autre personne est dans l'incapacité d'y accéder
pleinement. Néanmoins, on constate dans la pratique de la communauté
des mathématiciens que la validation d'un jugement est une notion
partageable: dans de très nombreux cas, les mathématiciens se mettent
d'accord pour reconnaitre qu'une inscription de démonstration permet
de valider un jugement, ce qui signifie qu'un travail d'exploration de
ce document génère une expérience de démonstration du jugement chez un
nombre suffisant de mathématiciens concernés.

Avec Wittgenstein, nous pouvons estimer que cette propriété ne relève
pas d'une explication philosophique mais d'un arrière-fond de
pratiques et de formes de vie \citep{Wittgenstein2004}.  À l'origine de
ce contexte de pratiques partagées, il y a l'enseignement de cette
science, dans lequel on apprend collectivement à éprouver l'évidence de
tel ou tel énoncé. En ce sens, l'apprentissage des mathématiques est
aussi un apprentissage collectif de la notion de démontrabilité et du
consensus qui en découle. Finalement, faire des mathématiques, c'est
être implicitement convaincu que la validité d'un raisonnement
mathématique est une notion qui a un sens collectif. D'ailleurs, il
est fréquent que la construction d'une démonstration soit l'œuvre non
pas d'une seule personne mais de tout un groupe de
protagonistes\footnote{Dans ce dernier cas, on peut estimer que le
  partage de l'expérience de l'évidence va au delà de la simple
  juxtaposition d'expériences indépendantes. En effet, l'intensité de
  la communication entre les interlocuteurs dans ses aspects verbaux
  et non verbaux (intonations, gestes, regards complices, sons divers,
  etc.)\ peut réussir à percer quelque peu l'extériorité de ces
  expériences en instaurant une forme de partage un peu plus
  conséquent. Cette situation est en partie analogue à celle qui
  aboutit au partage du monde physique qui nous entoure par
  l'intermédiaire de nos expériences perceptives, partage défendu
  fermement par James \citep{James2007} pour qui l'empirisme radical
  n'est pas un solipsisme.}.

Lorsqu'un jugement est partagé par un nombre significatif d'acteurs
compétents sur le sujet, la communauté des mathématiciens accorde sa
confiance à ce jugement. Cela implique que ce résultat devient une
sorte de «fait mathématique» sur lequel tout mathématicien peut
s'appuyer sans être dans l'obligation d'en éprouver lui-même
l'évidence.  Cette forme d'acceptation communautaire d'un nouveau
jugement est rendue nécessaire par l'étendue et la vitesse actuelle de
développement du champ mathématique. Sans doute que dans l'Antiquité,
même dans sa période tardive, l'extension du domaine mathématique
permettait à tout amateur sérieux d'en maitriser lui-même l'étendue
sans avoir à s'en remettre à l'avis de ses pairs.

Cette manière de procéder peut sembler être une forme de recours à un
oracle, mais elle en diffère parce que d'une part il y a eu partage de
la démonstration par un panel d'experts de confiance, et que d'autre part il
est toujours possible à quiconque de se reporter à la
démonstration. Ce dernier argument n'est pas tout à fait réaliste du
fait de la grande spécialisation des mathématiques actuelles et de
l'extrême complexité de certaines démonstrations. Quoi qu'il en soit,
le processus de validation d'un nouveau jugement est organisé en deux
niveaux se succédant dans le temps: le niveau individuel ou
quasi-individuel lorsqu'un ou quelques acteurs font l'expérience
initiale de l'évidence de ce jugement et en communiquent une
inscription de démonstration, puis le niveau collectif lorsqu'un
nombre significatif d'autres acteurs compétents partagent ce jugement
en approuvant l'inscription de démonstration. L'existence de ces deux
niveaux renforce considérablement la fiabilité des jugements
mathématiques.  Enfin, ce mécanisme de validation collectif en deux
étapes ne doit pas occulter le fait que, \emph{in fine}, c'est
l'expérience de démonstration qui constitue la pierre angulaire et la
spécificité de l'activité mathématique, même lorsque cette expérience est
déléguée à un groupe d'experts lors de la deuxième étape du processus
d'évaluation. En dernier ressort, c'est sur cette expérience que
repose toute la procédure.

\subsection{Le phénomène des erreurs \emph{post
    demonstrationem}} \label{erreurs}

Dans la vie courante, les raisonnements que nous faisons semblent
souvent bien fragiles: facilement remis en cause, ils sont la plupart
du temps l'occasion de querelles sans issues. La contemplation des
divers débats publics dans nos démocraties fatiguées en fournit un
exemple saisissant. \emph{A contrario}, les raisonnements
mathématiques ont la réputation d'être d'une solidité infaillible et
d'imposer de fait un consensus. C'est en partie justifié par la nature
même de la validation d'un jugement mathématique en deux étapes,
explicitées dans la partie~\ref{FormeDem} précédente, qui constituent
un réquisit d'une très grande exigence excluant en principe la
contestation. Néanmoins, l'idée d'un discours parfait tel que l'ont
pensé les philosophes platoniciens de l'Antiquité ou les laudateurs
contemporains de l'essence formelle des mathématiques, idée excluant
par essence toute faute, résulte d'un glissement abusif. La
manifestation la plus claire de cet abus de réputation réside dans le
phénomène des erreurs \emph{post demonstrationem}. Le caractère irrépressible de ce phénomène a été remarqué par Martin-Löf dans son
étude \citep{MartinLof87}. La terminologie «erreur \emph{post
  demonstrationem}» est introduite par les auteurs du présent
texte. Cette dernière ne doit pas être confondue avec la faute commise
par l'apprenti mathématicien qui ne sait pas encore articuler un
raisonnement mathématique. Le terme d'erreur \emph{post demonstrationem}
est employé lorsque: 1) à un moment donné, un jugement est validé par
la procédure décrite précédemment; 2) à un moment ultérieur,
éventuellement longtemps après, une erreur est découverte dans cette
validation. La manifestation de cette erreur est la découverte d'un
argument nouveau dont l'effet concret est de rendre impossible
l'expérience de l'évidence de ce jugement. Bien que ne mettant pas en
péril l'édifice des mathématiques, la fréquence d'apparition de ces
erreurs est cependant non négligeable. D'ailleurs, tout mathématicien
sait très bien qu'il lui arrive de trouver des erreurs dans des
raisonnements préalablement validés, que ces raisonnements proviennent
de lui-même ou d'une autre personne, qu'ils soient déjà publiés ou non
dans la littérature scientifique. Habituellement, la présence de ce
type d'erreur est analysée comme la simple marque de l'imperfection ou
de la finitude de l'esprit humain. Il en découle que les conséquences
philosophiques du phénomène des erreurs \emph{post demonstrationem} ne
sont pas correctement perçues.

Il est vrai que l'on est face à un paradoxe. D'une part, le processus
de la démonstration --~expérience de l'évidence, inscription guidant
cette expérience, partage de l'évidence guidée par l'inscription~--
est agencé, admirablement agencé peut-on ajouter, pour exclure les
erreurs: mener à son terme ce processus est la méthode, initiée en
grande partie par les mathématiciens grecs, qui apporte la garantie
que le raisonnement ainsi produit est valide, exempt de toute
faute. Mais d'autre part, la nature même de cette procédure ne permet
pas d'exclure définitivement la possibilité qu'une erreur puisse être
détectée par la suite. La raison en est qu'\emph{il n'y a aucun moyen
  extérieur à l'expérience de démonstration d'un jugement qui puisse
  permettre de s'assurer définitivement de sa validité}. \emph{La
  fiabilité d'une démonstration se fonde uniquement sur une forme
  d'expérience humaine, celle de l'évidence d'un jugement.} C'est déjà
beaucoup, surtout si l'inscription qui appuie cette expérience est
reconnue comme valide par une large communauté. Mais le raisonnement
mathématique ne peut pas prétendre atteindre une perfection sans
faille, une forme de validité éternelle, acquise pour la nuit des
temps, car rien dans ce qui constitue la pratique des mathématiques ne
peut justifier cette prétention. Une simple expérience humaine n'étant
pas taillée pour atteindre l'absolue perfection, il faut se résoudre à
reconnaitre une irréductible part de fragilité aux résultats
mathématiques.

Ainsi, l'assurance d'une garantie absolue du bien-fondé d'un jugement
transcendant l'expérience de démonstration est un mythe et un
non-sens. La technique de la démonstration d'un jugement est la
méthode mise au point par les hommes pour attester de la validité de
ce jugement; faute d'être acquise à tout jamais, cette validité est
assurée dans le temps de l'expérience humaine, c'est-à-dire tant qu'il
y aura des hommes qui exploreront la démonstration de ce jugement pour
en éprouver l'évidence et qui confronteront ce jugement avec les
nouveaux acquis mathématiques avec lesquels il pourrait être en
relation. D'ailleurs, il est notable que les mathématiciens préfèrent
généralement revisiter les démonstrations des jugements qu'ils doivent
utiliser dans leur propre travail plutôt que de les admettre
directement\footnote{Bien entendu, il ne s'agit que d'une préférence,
  et un mathématicien n'a pas toujours la possibilité de retraverser
  tous les outils qu'il utilise. C'est le cas pour certains théorèmes
  de nos mathématiques actuelles dont la démonstration complète
  nécessite plusieurs milliers de pages comme par exemple le «théorème
  énorme» qui classifie les groupes simples finis.}.

D'une certaine manière, c'est la croyance que la technique de la
démonstration n'est qu'un moyen pour accéder à une instance extérieure
à cette expérience qui fonde l'idée de l'incontestabilité du
raisonnement mathématique. Pour les philosophes platoniciens, cette
instance était le monde transcendant des idées pures. Pour les
modernes d'aujourd'hui, cette instance est plus volontiers le monde
des objets formels. D'où cet ultime refuge pour la thèse de
l'infaillibilité du discours mathématique: un raisonnement formalisé,
c'est-à-dire complètement coulé dans le moule syntaxique de l'une des
théories formalisant la pratique des mathématiques, peut sembler par
essence complètement imparable. Si on admet cette thèse, on est alors
tenté de voir à nouveau tout raisonnement usuel comme une
approximation nécessairement imparfaite d'un objet formel parfait
(voir la citation de Bourbaki page~\pageref{citation-bourbaki}), ce
qui expliquerait à la fois la présence d'erreurs cachées dans la
pratique ordinaire des mathématiques et une stratégie pour les éviter,
à savoir celle de s'approcher autant que faire se peut de la
formalisation complète. Mais qu'est-ce qu'une démonstration formelle?
Comme cela a déjà été noté, une configuration de signes sur un support
quelconque ne peut prendre le statut de démonstration formelle que
dans le cas où elle est identifiée comme telle dans le cadre d'une
expérience d'un sujet mathématicien. Quelle que soit la nature de
cette expérience, il ne s'agit que d'une simple expérience humaine,
inapte à fonder une infaillibilité sans restriction.  En conséquence,
la stratégie consistant à assoir la perfection absolue d'un
raisonnement mathématique aboutissant à un certain jugement en le
formalisant complètement ne peut pas fonctionner: à supposer que l'on
soit capable de procéder à cette formalisation, cela ne ferait que
déplacer la possible intrusion future d'erreurs dans l'expérience qui
assure que l'entité syntaxique obtenue possède la qualité de
démonstration formelle\footnote{Ce qui ne veut pas dire que le travail
  sur les démonstrations formelles ne présente pas d'intérêt
  scientifique!  C'est seulement une certaine utilisation de ces
  développements formels à des fins philosophiques qui est critiquée
  ici.}.

Finalement, il n'existe pas de méthode permettant de s'assurer
définitivement de la validité absolue d'une démonstration. D'ailleurs,
s'il en existait une, il suffirait de l'appliquer à chaque
démonstration produite pour être certain qu'elle ne sera jamais mise
en défaut par le surgissement d'une erreur. Néanmoins, il ne faut pas
en déduire que le raisonnement mathématique est miné par des erreurs
potentielles. Ce que révèle l'analyse précédente, c'est que le
mécanisme de la preuve ne permet pas d'exclure à tout jamais la
possibilité d'une erreur, ce qui ne veut pas dire qu'il y aura dans le
futur --~certainement ou avec une grande probabilité~-- apparition
d'erreurs. De fait, en l'absence d'une instance transcendante
permettant d'évaluer dans l'absolu un jugement donné, une expérience
de démonstration de ce jugement offre une excellente garantie sur sa
validité, à la fois extrêmement solide et limpide sur son
fondement\footnote{Cette garantie est quelque peu mise à mal par la
  course à la publication qui anime la recherche actuelle. Même en
  mathématiques, cette fièvre productiviste pousse à ne pas trop se
  poser de question sur la fiabilité des outils utilisés et hélas
  aussi sur la qualité des preuves. D'où des erreurs fréquentes qui
  sont relevées par les relecteurs des articles, du moins lorsque ces
  derniers acceptent de dégager le temps nécessaire à cette tâche, au
  détriment de leurs propres publications\dots}.

\subsection{Deux remarques} \label{DeuxRem}
 
\paragraph{1.}  Jusqu'à ce point, l'analyse présentée auparavant n'a
nullement pris en compte les fortes différences, pratiques et
théoriques, qui en gros scindent la communauté des mathématiciens en
deux écoles, d'une part les mathématiques classiques, et d'autre part
les mathématiques constructives\footnote{En termes d'effectif, ces
  deux écoles sont inégales puisque seulement une minorité de
  mathématiciens peuvent être qualifiés de constructivistes. Cette
  disparité n'empêche pas le courant constructiviste d'avoir à nouveau
  une grande importance, en particulier du fait de ses liens forts,
  théoriques et pratiques, avec l'informatique et la logique.}. En
fait, la présentation d'une conception non formaliste de la démonstration
présentée ci-haut est censée s'appliquer de manière égale aux
mathématiques développées dans ces deux groupes. Cela peut paraitre
étonnant à un lecteur connaisseur de l'histoire des mathématiques
depuis la fin du 19\ieme~siècle, qui sait parfaitement que ces deux
courants s'opposent en particulier sur ce que doit être une
démonstration. De plus, ce même lecteur peut aussi remarquer que la
définition de la notion d'expérience de démonstration donnée dans la
partie~\ref{PreuveExpe} semble relativement proche des conceptions
philosophiques développées dans l'intuitionnisme, comme l'indique la
note~\ref{IntuiML} faisant référence à Martin-Löf. Cependant, les
auteurs du présent texte défendent que, au delà de l'opposition
philosophique sur la querelle des fondements, c'est la forme
paradigmatique de la notion de démonstration, telle qu'elle est
pratiquée depuis les Grecs, qui est bien captée par la définition
choisie ici. Tout mathématicien, classique ou constructiviste, admet
la validité d'une démonstration lorsque, au terme d'un travail
d'exploration et de maturation, il est saisi par l'évidence de la
conclusion. Jusqu'à maintenant, cette expérience de la démonstration
est constitutive de l'essence de l'activité mathématique, et elle
transcende les controverses. C'est dans le choix des moyens utilisés
pour édifier la démonstration que l'opposition est
pertinente\footnote{L'opposition est aussi totale en ce qui concerne
  la notion de vérité (voir par exemple la présentation donnée dans
  \citep{Bourdeau2004}). Cette dissension majeure ne semble pas
  affecter le fait qu'une démonstration est vécue comme l'expérience
  de l'évidence de sa conclusion.}. Par exemple, l'utilisation de la règle du tiers exclu est jugée légitime par un mathématicien
classique au sens où elle lui semble marquée par l'évidence, alors que
cette utilisation est rejetée par un mathématicien constructiviste
pour qui ce même moyen rompt le cheminement de l'évidence\footnote{À
  ce propos et en dehors de tout engagement partisan, les auteurs du
  présent texte ont pour leur part le sentiment que, en ce qui
  concerne le cheminement de l'évidence qui constitue le cœur de
  l'expérience de démonstration, les mathématiques constructives
  constituent un cadre plus favorable à la force de cette
  expérience.}.

\paragraph{2.}  Une autre question non encore abordée dans cette
discussion est relative à l'utilisation de calculs sur machine dans
une démonstration. En effet, de manière analogue à ce qui se passe
dans toutes les autres sciences, les progrès technologiques et
théoriques font apparaitre des dispositifs de calcul susceptibles
d'être utilisés dans la vérification ou la construction de
démonstrations. C'est le cas des assistants de preuve comme le
logiciel \emph{Coq} \citep{coqmanual} qui permettent de vérifier, de
compléter ou même de construire des démonstrations. Un exemple
maintenant classique de la réussite de cette démarche est le théorème
des quatre couleurs dont la démonstration complètement certifiée a été
obtenue en 2004 \citep{gonthier05} et dont on ne possède aucune démonstration pouvant se
passer de calculs sur machine.

Pour illustrer l'impact de l'usage des ordinateurs par un exemple
simple, imaginons d'abord que la seule démonstration que nous aurions
de ce que la somme des nombres entiers de~1 à 2~millions fait
2~billions et 1~million était l'inscription du calcul de cette somme
sur une machine; or nous avons une démonstration, qui est
que cette somme est aussi celle des nombres entiers de 2~millions à~1,
et qu'en additionnant terme à terme les termes de ces deux sommes, on
obtient 2~millions de fois 2~millions et~1, qu’il reste à diviser par deux. L'évidence de cette
démonstration nous parait d'une autre qualité que l'inscription du
calcul.

La question est de savoir en quoi la notion de démonstration présentée
auparavant pourrait être modifiée par ces nouvelles pratiques. Par
exemple, considérons une démonstration dont une partie n'est validée
que par un calcul sur machine; est-ce qu'elle peut générer une démonstration au
sens de l'expérience de l'évidence de son résultat?  À première vue,
il semble que non parce que le cheminement de l'évidence guidé par la
démonstration est \emph{a priori} rompu: la partie de la démonstration
prise en charge par la machine se présente comme une béance
à ce niveau, du moins pour celui qui n'est pas capable de comprendre
ce que fait le programme informatique. Notons que l'intrusion de la
matérialité dans l'efficience des mathématiques n'est pas nouvelle et
elle est même structurelle d'après la contribution de Netz sur
l'émergence des mathématiques dans le monde grec; ce qui semble
nouveau avec la mécanisation des démonstrations est que le dispositif
mis en œuvre certes peut apporter une réponse, mais que cette dernière
n'est pas de nature à favoriser la démonstration au sens de
l'expérience de l'évidence de sa conclusion.

Néanmoins, un mathématicien peut être en capacité de saisir le fonctionnement
de la machine, d'en démontrer la correction, pour finalement estimer que, tant bien que mal,
l'évidence poursuit son chemin dans la partie traitée
machinellement. Encore faut-il avoir une confiance absolue en le
dispositif électronique piloté par le programme qui effectue la tâche
demandée. Il faudrait une confiance de l'ordre de celle que l'on
accorde à un résultat mathématique validé par une démonstration,
c'est-à-dire, une confiance qui exclut totalement le doute du fait que
l'on est pénétré de son évidence par le contrôle que l'on a pu exercer
sur l'intégralité du processus qui le valide.  On voit bien
que cela ne peut pas fonctionner ainsi. Dans le cas considéré, le
détail du fonctionnement intime du processus machinel qui valide le
résultat n'est justement pas accessible (nous négligeons ici une
défaillance toujours possible au niveau matériel de la mise en œuvre
par une machine concrète). Comme indiqué dans la la
note~\ref{NoteOracle}, on est dans le cas de l'utilisation explicite
d'un oracle selon la terminologie utilisée à la fin de la
partie~\ref{PreuveExpe}.

Force est de constater que les mathématiques tributaires de telles
démarches informatiques se rapprochent des autres sciences comme la
physique qui dépendent de manière plus ou moins essentielle de
dispositifs matériels expérimentaux. Mais alors il faut être conscient
que ces mathématiques nouvelles sont, dans leur méthode même, d'une
nature différente du procédé de la démonstration instauré par les
Grecs, qui depuis ce moment inaugural a constitué l'originalité et la
marque distinctive de cette discipline.  Il est possible que ce soit
une nouvelle science qui émerge de cette manière, à la frontière des
mathématiques, de la logique et de l'informatique: une science dont
les objets sont les démonstrations formelles, et dont la démarche est
fortement dépendante d'outils informatiques puissants. La confiance
accordée aux résultats fournis par ces démarches peut être très forte
car le travail scientifique concerné est souvent de qualité
incontestable. Par exemple, la certification \emph{Coq}, label donné à
un résultat vérifié par l'assistant de preuve \emph{Coq}, est
considérée comme une forme très haute de validité\footnote{Plus
  précisément, la certification a lieu sur la base d'une inscription
  de démonstration spécialement adaptée à l'assistant de preuve,
  élaborée à partir d'une démonstration dont un mathématicien a fait
  l'expérience et qu'il a rédigée. Alors pourquoi certifier? Trois
  raisons nous viennent à l'esprit: les erreurs \emph{post
    demonstrationem} (c'est la motivation initiale de Voevodsky, qui a
  été très affecté par de telles erreurs comme évoqué dans la
  note~\ref{note:voevodsky}); convaincre les autres mathématiciens de
  la correction de la démonstration et en particulier de l'absence d'oubli
  d'un cas particulier parmi un grand nombre d'autres (théorème de
  Hales de la conjecture de Kepler, théorème des quatre couleurs);
  développer la certification en tant que discipline scientifique
  (théorème de Feit-Thompson; on peut noter qu'au cours de la
  certification de ce théorème, quelques erreurs de la démonstration
  originelle ont été corrigées, mais qu'elles étaient superficielles).}, à tel point que si un résultat obtenu par
\emph{Coq} entrait en contradiction avec un jugement validé par une
démonstration au sens développé précédemment de l'évidence vécue, il
ne va pas de soi que la confiance serait prioritairement accordée au
deuxième. Autrement dit, sur le plan de la force de la confiance, un
résultat certifié par \emph{Coq} n'a rien à envier à un résultat
certifié par une démonstration.  Mais ce n'est pas la force en elle-même de
la confiance accordée à un résultat qui est la marque distinctive de
l'essence de la démarche mathématique. C'est la nature du processus
qui légitime le résultat, à savoir la démonstration comme expérience
de son évidence qui est spécifique aux mathématiques depuis les
Grecs. Cette expérience crée chez celui qui l'éprouve une forme de
rapport particulier avec le résultat validé, rapport qui est
incompatible avec la délégation d'une partie de la démonstration à un
oracle.

\section{Conclusion: vers un mode d'existence pour les mathématiques}
 
Cette partie terminale a pour objet principal d'introduire une forme
de mode d'existence propre aux êtres mathématiques respectant le patron de l'\textsc{eme} et de discuter le
placement de ce mode dans celui de la référence
[\textsc{ref}]. Préalablement, quelques précisions doivent être
apportées relativement au mode [\textsc{ref}] et à la notion générale
de mode d'existence au sens de l'\textsc{eme}.
  
\subsection{Quelques précisions sur les modes d'existence en général
  et sur le mode [\textsc{ref}] en particulier} \label{PrecisREF}
 
Il n'est guère utile pour un lecteur ne connaissant pas le sujet de
donner une description de la notion de mode d'existence selon
l'\textsc{eme} en égrenant les termes du métalangage qui a été forgé
pour fixer un cadre général à cette notion. Il est certainement
préférable d'évoquer préalablement un cas concret. C'est ce qui va
être fait avec la présentation succincte de la manière dont le travail
d'une mission scientifique de terrain est analysé par l'anthropologie
des sciences. Le type de pratique scientifique correspondant à cette
mission semble de prime abord aussi éloigné que possible de la
pratique des mathématiques, ce qui devrait permettre par la suite de
mieux faire apparaitre les analogies et différences entre les
modalités d'existence des êtres sous-jacents à ces pratiques.

L'exemple choisi est celui de la lumineuse étude\footnote{Cette étude
  est un splendide exercice de philosophie appliquée. La description
  qui en est donnée dans le présent article est d'une brièveté qui ne
  lui rend pas justice.} développée dans
\citep{Latour2007}. Ce texte analyse le travail d'une mission
scientifique dont l'objet est de comprendre la dynamique de la
transition forêt-savane dans une région reculée du Brésil proche de
Boa Vista. L'étude montre l'élaboration concrète d'une chaine de
médiations --~maintenant appelée chaine de référence~-- depuis le
choix de la parcelle choisie pour mener l'investigation, la
triangulation de cette parcelle à l'aide d'un «pédofil\footnote{C'est
  le petit nom donné à un instrument permettant de mesurer les
  distances dans un terrain naturel encombré et non aplani.}», les
prélèvements d'échantillons de terre à différentes profondeurs et leur
placement dans une pédocomparateur\footnote{Dispositif permettant de
  ranger des échantillons de terrain selon une disposition planaire en
  rapport avec la localisation des prélèvements.}, les traitements de
ces échantillons permettant d'extraire du pédocomparateur des tableaux
numériques, enfin l'élaboration d'un diagramme final tenant sur une
simple feuille de papier et résumant à lui tout seul l'investigation,
la conclusion de cette dernière sur la structure du sol au voisinage
de la lisière et la formulation de l'hypothèse d'une probable avancée
dans le temps de la forêt sur la savane grâce au laborieux travail des
vers de terre. Cette chaine est donc constituée par la succession
d'inscriptions scripto-visuelles, chacune étant plus abstraite et
moins matérielle que la précédente. Le choix judicieux de ces formes,
le protocole très précis avec lequel ce choix est fait ainsi que
l'opération qui fait passer d'une inscription à la suivante, tout cela
permet à l'information de circuler dans les deux sens le long de la
chaine, entre les points extrêmes que sont la parcelle initiale et le
diagramme final. L'évolution des inscriptions le long de la chaine va
vers la simplicité perceptive en passant graduellement de
l'indéchiffrable fouillis du terrain initialement choisi au limpide
diagramme terminal. En effet, la fonction de la chaine est d'arriver à
une inscription finale qui soit suffisamment claire, simple et
significative pour emporter la conviction des chercheurs
concernés. Finalement, ce que l'on pourrait nommer de manière quelque
peu pédante «l'ontologie objective de la nature du sol
transversalement à la lisière» est ce que donne à voir le diagramme
final en gardant soigneusement en mémoire la circulation de la
référence tout le long de la chaine.
 
Cette description met concrètement en évidence les traits principaux
des chaines de référence. Selon l'\textsc{eme}, ces chaines
d'inscriptions sont au cœur du mode d'existence [\textsc{ref}], au sens
où quelque chose existe objectivement lorsque l'on peut mettre en
place une chaine de référence qui établit cette chose \emph{via} la
cascade d'inscriptions constituée par cette chaine. La dissemblance,
l'éloignement entre chaque inscription et l'inscription suivante, est appelée le
\emph{hiatus} propre au mode [\textsc{ref}]. Ce hiatus peut être perçu
au niveau très général de la seule dissemblance de deux inscriptions
en tant que simples entités scripto-visuelles (dissemblance que l'on
pourrait interpréter comme des différences entre les parties
textuelles ou graphiques de chacune des inscriptions). De manière plus
spécifique, ce hiatus peut aussi signifier le manque \emph{a priori}
de conviction que la seconde inscription est conçue de manière à
porter une information que l'on peut rapporter à la première, cette
conviction étant une expérience recherchée par les acteurs
scientifiques qui élaborent la chaine ou en prennent connaissance.

La \emph{condition de félicité} qui permet de lisser une telle chaine
est que l'opération qui fait passer d'une inscription à la suivante
soit conçue de manière à préserver une certaine constance\footnote{On
  peut remarquer l'analogie avec une méthode utilisée pour montrer
  qu'un algorithme itératif est correct, c'est-à-dire qu'il fait bien
  ce qu'on attend de lui. Cette méthode consiste à chercher un
  invariant de boucle, c'est-à-dire un prédicat qui reste constant au
  cours des itérations.}, de telle sorte que l'information lisible sur
la deuxième inscription remonte à une information portée par l'inscription
précédente. La dernière inscription d'une chaine de référence permet
en principe de conclure mais elle n'a de valeur que par l'ensemble de
la chaine qui, du fait des conditions de félicité, permet de faire
transiter cette information en remontant pas à pas la chaine de
manière à aboutir au phénomène initial sur lequel porte l'étude.

Après la présentation d'un exemple concret et de certains éléments du
mode [\textsc{ref}] correspondant, il est temps d'énoncer quelques
généralités significatives sur les modes d'existence au sens de
l'\textsc{eme}. Un mode d'existence désigne un collectif d'entités
partageant une même manière d'être et de s'animer, et dont la présence
active peut être décelée dans notre modernité, du moins lorsqu'on
abandonne le dogme selon lequel il y a une seule forme d'existence qui
est celle de la vérité objective à laquelle la Science et la Raison
donnent un accès direct et sans perte. L'une des idées-phares de
l'\textsc{eme} est que, bien qu'une multiplicité de modes d'existence
puissent être distingués, ces derniers partagent des traits communs
qui permettent de les identifier et de les classer. Ce sont
principalement:
\begin{itemize}
\item un certain type de \emph{trajectoire} qui est le support
  principal de l'existence dans le mode considéré --~dans le cas de
  [\textsc{ref}], ce sont les chaines de référence;
\item un certain type de \emph{hiatus} qui sépare deux étapes
  successives quelconques d'une trajectoire --~dans le cas de
  [\textsc{ref}], c'est l'éloignement, la dissemblance entre deux
  inscriptions successives d'une chaine de référence;
\item un certain type de \emph{conditions de félicité et d'infélicité}
  qui sont les conditions de vérité ou de fausseté propres au mode
  considéré --~dans le cas de [\textsc{ref}], étant donnée une chaine
  de référence et une information qui émerge à la dernière étape de la
  chaine, la condition de félicité énonce que cette information est
  pertinente si elle peut être rapportée d'étape en étape en remontant
  toute la chaine de référence.
\end{itemize}
 
\subsection{Un quasi-mode d'existence pour les
  mathématiques} \label{quasimode}
 
Le propos de cette partie est de décrire le fonctionnement, la manière
d'être et de se présenter à nous, des êtres qui se manifestent quand
un mathématicien (ou un groupe de mathématiciens) met au point un
développement mathématique ou parcourt un tel développement. Comme
cette description peut se faire en suivant le patron des modes
d'existence au sens de l'\textsc{eme}, il semble naturel d'en déduire
que les entités mathématiques relèvent d'un tel mode. En toute
rigueur, il serait préférable de parler pour l'instant d'un quasi-mode
car, à supposer que l'argumentaire qui suit soit suffisant pour
établir une forme de mode d'existence acceptable, il restera à le
placer dans l'architecture générale proposée par l'\textsc{eme} selon
l'une des options suivantes: comme relevant purement et simplement du
mode de la référence objective [\textsc{ref}], ou bien comme sous-mode
de ce même [\textsc{ref}] mais doté d'une forte spécificité, ou bien
comme nouveau mode à part entière, ou bien comme autre chose
encore\dots{} L'idée qui va être développée est que ce quasi-mode
d'existence est déjà presque présent, en pointillé, dans la
description empiriste du raisonnement mathématique donnée dans la
partie~\ref{NonForm}.
 
Préalablement, il convient de revenir sur les inscriptions
scripto-visuelles qui interviennent dans la pratique des
mathématiques.  Dans la partie~\ref{FormeDem}, avec la notion
d'inscription de démonstration (d'un jugement), nous avons introduit
une telle inscription en insistant sur la différence essentielle avec
l'expérience de démonstration (du même jugement). Il est maintenant
utile de considérer la notion \emph{a priori} banale d'\emph{énoncé
  mathématique} (expression simplifiée en \emph{énoncé} dans la
suite): habituellement, c'est un texte qui exprime un contenu
mathématique, que ce dernier soit établi --~ce qui en fait un
jugement~-- ou ne le soit pas. Élargissons cette
notion de manière à inclure éventuellement d'autres éléments afin d'en
faire un tout relativement complet sur le plan de l'expression
habituelle des mathématiques. Par exemple, l'ajout peut être une
partie graphique, un diagramme, une image, etc. Évidemment, tous ces
éléments extra-textuels doivent être comptés à égalité d'importance
avec la partie strictement textuelle.  Pour prendre un exemple
élémentaire et bien connu, le mathématicien qui évoque le théorème de
Thalès met en place immédiatement, sur un support matériel et/ou dans
une image mentale, un texte et un diagramme avec lettres: le tout
décrit ce résultat et en constitue l'accrochage matériel et
symbolique. Ce cadre scripto-visuel --~cette forme~-- est pratiquement
indispensable à tout travail sur le résultat considéré. Dans la suite,
l'expression d'énoncé est prise au sens élargi que nous venons de
préciser. Après ce préalable, le mode d'existence pour les
mathématiques peut être abordé en précisant son type de hiatus, de
trajectoire et de condition de félicité.
 
Le hiatus du quasi-mode d'existence pour les mathématiques est la
discontinuité manifestée par l'éloignement entre un énoncé donné
initialement et un énoncé final. Cet éloignement peut être compris
comme la simple dissemblance entre les deux inscriptions
scripto-visuelles en jeu; de manière plus spécifique, il peut aussi
être interprété comme l'absence d'expérience de démonstration qui permettrait de passer du
premier au second.

Le saut, la transition, la transformation reliant ces énoncés est une
inscription de démonstration qui est présentée comme un guide susceptible de faire transiter
l'évidence du premier au second. Pour important qu'il soit, ce guide
ne certifie pas que l'expérience du cheminement de l'évidence sera au
rendez-vous de celui qui le parcourt.

La condition de félicité de ce saut est que cette démonstration soit
considérée comme valide, c'est-à-dire qu'elle génère l'évidence de
l'énoncé final chez celui qui éprouve cette «passe» \citep[«Livre»,
colonne «Vocabulaire», entrée «Passe»]{EME} particulière, c'est-à-dire l'expérience du cheminement de l'évidence du premier au second
énoncé. Lorsque
cette condition est réalisée, l'énoncé final par delà le hiatus
devient un jugement, sous réserve que l'énoncé source en deçà soit
lui-même un jugement, et le saut s'identifie à une
inscription de démonstration de ce
jugement à partir du jugement source. Autrement dit, cette condition
de félicité est satisfaite lorsque l'évidence peut circuler jusqu'à
l'énoncé final en suivant le guide fourni par la démonstration.

Une fois définis les éléments précédents, on considère les
trajectoires obtenues en enchainant hiatus et sauts, trajectoires qui
seront nommées au moins provisoirement \emph{chaines
  démonstratives}. Précisément, une telle chaine est donnée par:
\begin{itemize}
\item une succession d'inscriptions $I_1,I_2,\ldots,I_n$ qui sont des
  énoncés; la première de ces inscriptions~$I_1$ est supposée être un
  jugement;
\item une succession d'inscriptions $D_1,D_2,\ldots,D_{n-1}$ telles
  que pour chaque $k=1,2,\ldots,n-1$, l'inscription~$D_k$ soit une
  démonstration de~$I_{k+1}$ à partir de~$I_k$.
\end{itemize}
On peut représenter ceci sous la forme du graphique suivant:
\begin{equation}\label{Chaine1}
  I_1 \mapright{D_1} I_2 \mapright{D_2} I_3 \mapright{D_3}
  \cdots \mapright{} I_{n-1} \mapright{D_{n-1}} I_{n}.\tag{*}
\end{equation}
 
De fait, cette notion de chaine démonstrative n'est qu'une certaine
mise en forme, en termes de chaines d'inscriptions, de la notion
générale de démonstration. Cette manière de présenter la notion de
démonstration a pour objet d'en donner une présentation proche de
celle d'une chaine de référence sans en déformer le contenu. Avec les
notations introduites dans la définition~(\ref{Chaine1}), la chaine
démonstrative peut être vue comme la démonstration de~$I_{n}$ à partir
de~$I_{1}$. D'un point de vue purement formel, cette démonstration
globale pourrait remplacer la chaine toute entière en réduisant cette
dernière à un seul saut
\[
  I_1 \mapright{D} I_n\text.
\]
Dans la pratique, il en va autrement, car, dans le but de faire
circuler plus facilement l'évidence, le mathématicien sait qu'il est
préférable de scinder la démonstration globale en une succession de
démonstrations intermédiaires, chacune étant plus facile à élaborer et
plus à même de susciter le transport de l'évidence\footnote{En ce
  sens, l'argumentation mathématique n'est pas aussi éloignée de la
  rhétorique que l'on se plait à le dire.}. Cela justifie pleinement
la notion de chaine démonstrative, qui plus est en exhibant une
propriété importante: pour que ces chaines deviennent opérationnelles,
il est préférable que les énoncés successifs soient suffisamment
«proches» au sens où il est relativement aisé de faire transiter
l'évidence entre eux. À propos de la structure essentiellement
linéaire d'une telle chaine --~certainement réductrice par rapport à
la complexité combinatoire que peut prendre en général une
démonstration~-- il est bon de rappeler la mise en garde exprimée dans
la note~\ref{StructDem} (et les mises en garde de Gentzen
\cite[pages~19-20]{Gentzen1935} lui-même).
 
La condition de félicité définie précédemment pour un saut isolé se
généralise dans les mêmes termes à une chaine toute entière: une
chaine démonstrative satisfait la condition de félicité lorsque chacun
de ses sauts constitutifs la satisfait, autrement dit lorsque
l'évidence peut cheminer le long de la chaine en validant comme
jugement chacun des énoncés successifs. Au contraire, la condition
d'infélicité est satisfaite lorsque le cheminement de l'évidence est
interrompu en un point de cette chaine.

Finalement, le quasi-mode d'existence des êtres mathématiques est
porté par ces chaines démonstratives lorsqu'elles sont reconnues
satisfaire la condition de félicité. Un être mathématique
(c'est-à-dire un objet ou une propriété mathématique) est amené à
l'existence lorsqu'une chaine démonstrative satisfaisant la condition
de félicité énonce, montre cette existence. C'est ainsi que tout
jugement, habituellement conçu comme un acquis de connaissance, est
également interprétable comme une extension de l'ontologie: telle
entité (objet ou propriété) est dotée d'une existence avérée par la
chaine démonstrative correspondante. Par exemple, une chaine
démonstrative associée au théorème de Pythagore amène à l'existence,
pour un triangle rectangle, la propriété d'égalité entre les carrés
construits sur les côtés adjacents à l'angle droit et le carré
construit sur l'hypoténuse.

Dans quel lieu ces êtres vivent-ils? L'ambigüité ontologique des
mathématiques relevée par Netz fait qu'aucun lieu n'est attribué
d'emblée aux objets mathématiques dans ce que nous nommons naïvement
la réalité extérieure. Depuis les mathématiques grecques, la réponse
platonicienne à cette difficulté a été de postuler que ces objets se
tiennent figés dans l'immuable éternité d'une sorte de monde
suprasensible des idées pures. Le quasi-mode d'existence des êtres
mathématiques présenté dans cette partie~\ref{quasimode} permet
avantageusement d'éviter ce recours à la transcendance tout en restant
au plus près de la pratique des acteurs de cette sciences.  Il reste
que, si on le souhaite, on peut dire que les entités mathématiques
peuplent des mondes qui leur sont propres tout en restant attachées
aux dispositifs scripto-visuels qui les supportent. Être
mathématicien, c'est peut-être d'abord avoir la capacité de
reconnaitre ces mondes, d'y accéder, de les voir s'animer à partir du
seul travail sur des inscriptions appropriées. Du fait du caractère
spécifique et généralement difficile à mettre en œuvre du processus
qui amène ces êtres à l'existence, il est naturel de les qualifier de
lointains.  Le seul moteur d'animation d'un monde d'entités
mathématiques est le cheminement de l'évidence le long des chaines
démonstratives correspondantes. Dans ce nouveau cadre, il ne faut pas
oublier la leçon de Netz et de Peirce: un tel monde ne devient
opérationnel que parce que certaines conditions techniques et
cognitives sont satisfaites, conditions qui permettent à la pratique
des mathématiques de se déployer principalement dans la dimension
diagrammatique (voir la partie~\ref{NetzPeirce}). Évidemment, ces
conditions accentuent encore plus l'éloignement d'un tel monde
relativement à celui des «affaires courantes».

Comme cela est énoncé dans l'\textsc{eme}, «\emph{chaque saut par-dessus une discontinuité représente un risque pris qui peut réussir
  ou rater}». Dans le cas qui nous intéresse, le risque est que l'inscription de
démonstration soit considérée comme insuffisante pour générer
l'expérience de l'évidence de l'énoncé final.  L'objet de la condition
de félicité est donc de surmonter ce risque. Cependant, le risque ne
peut pas être éliminé définitivement: c'est la leçon enseignée par le
phénomène des erreurs \emph{post demonstrationem}, c'est-à-dire la
nécessité d'envisager l'éventualité d'erreurs dans la démonstration du
jugement considéré. Traduite dans notre perspective ontologique, cette
leçon montre que l'existence d'un être mathématique ne peut pas avoir
la valeur d'une vérité définitive transcendant l'expérience du
cheminement de l'évidence éprouvé par un mathématicien. Cet être est
doté d'une existence certaine tant que le transport de l'évidence qui
le valide est reconnu comme tel, ce qui ne peut exclure la découverte
future d'une faille susceptible d'invalider cette existence. La raison
en est que ce cheminement est une expérience, celle de l'évidence du
jugement final, et non pas l'accès à une forme de réalité indépendante
des acteurs mathématiciens. Comme indiqué précédemment, on peut y voir
la raison de la manie qu'ont les mathématiciens de parcourir les
argumentaires menant à des êtres mathématiques pourtant déjà acquis:
ils testent le risque d'une erreur et font ainsi l'expérience de leur
existence. Le maintien dans l'existence --~la \emph{subsistance} selon
le vocabulaire général de l'\textsc{eme}~-- des êtres mathématiques
pourrait se faire au prix de cette vérification sans cesse renouvelée.

Dans la suite de ce travail, le quasi-mode précédent va être comparé
au mode de la référence objective tel qu'il a été présenté dans la
partie~\ref{PrecisREF}.  Cette comparaison permettra de discuter
quelle place pourrait être attribuée aux êtres mathématiques dans
l'\textsc{eme}.

\subsection{Éléments de comparaison entre [\textsc{ref}] et le
  quasi-mode d'existence des êtres mathématiques} \label{CompModes}
 
La comparaison entre le quasi-mode d'existence des êtres mathématiques
et celui de la référence objective se réduit principalement à la
comparaison des chaines et des conditions de félicité propres à
chacun. C'est ce qui va être fait dans la suite en mettant en lumière
des analogies et des différences entre
\begin{itemize}
\item d'une part la technique des chaines de référence pratiquée dans
  les sciences de terrain ou expérimentales rapidement présentées dans
  la partie~\ref{PrecisREF} et qui, pour éviter un risque de confusion
  avec une éventuelle notion plus large de chaine de référence, seront
  provisoirement appelées des \emph{chaines de référence objective},
\item d'autre part le travail de la preuve mis en œuvre en
  mathématiques et représenté par les \emph{chaines démonstratives}
  introduites dans la partie~\ref{quasimode} précédente.
\end{itemize}
 
La première analogie est que ces deux types de chaines se présentent
comme des trajectoires dont les nœuds --~les étapes successives~--
sont des inscriptions scripto-visuelles. Il s'agit d'une analogie
relative au cadre général des modes d'existence décrits dans
l'\textsc{eme} puisque seules les trajectoires du mode [\textsc{ref}]
présentent cette particularité. Pour autant, la nature de ces
inscriptions peut révéler une différence notable entre les deux types
de chaines. En effet, comme c'est le cas dans l'exemple présenté dans
la partie~\ref{PrecisREF}, les inscriptions d'une chaine de référence
objective peuvent inclure des données matérielles extraites du milieu
étudié, ce qui n'est pas le cas d'une chaine démonstrative. Cependant,
cette différence tend à se gommer sinon à disparaitre lorsque l'on
considère le déroulement d'une chaine de référence objective. En
effet, les inscriptions d'une chaine démonstrative sont des énoncés
qui utilisent les outils habituels de représentation d'entités
mathématiques par un mixte de textes et d'images. Or, comme cela est
relevé dans la partie~\ref{PrecisREF}, l'une des propriétés des
chaines de référence objective est l'utilisation croissante le long de
la chaine de ces mêmes outils. Précisément, la succession des formes
d'une chaine de référence objective est marquée par une double
dynamique: la diminution, allant le plus souvent jusqu'à la
disparition pure et simple, de données matérielles en rapport avec le
phénomène étudié, et l'augmentation d'entités de nature mathématique.
 
Une deuxième analogie significative apparait au niveau de la fonction
des chaines de référence et des chaines d'existence
mathématique. Dans les deux cas, ces dispositifs ont pour but de créer
une conviction à propos de la vérité d'un énoncé scientifique, cette
conviction étant atteinte à la fin de la construction de la chaine
dans les deux cas. Comme la suite va le montrer, cette ressemblance ne
va pas jusqu'à inclure la nature même de cette conviction.
 
Une troisième analogie concerne le hiatus se présentant dans chacun
des deux types de chaines, hiatus qui dans les deux cas est formulé de
manière identique comme la dissemblance et l'éloignement de deux
inscriptions successives. Comme précisé précédemment, cet éloignement
peut être interprété comme un manque \emph{a priori}, celui de
l'expérience en laquelle consiste la satisfaction de la condition de
félicité correspondante. L'analogie en question concerne donc la
formulation du hiatus, ce qui n'exclut pas une certaine dissemblance
découlant elle-même d'un contraste entre les deux types de conditions
de félicité.

Justement, une quatrième analogie peut être relevée en ce qui concerne
les conditions de félicité.  Dans chacune de ces chaines, chaque saut
d'une inscription à la suivante doit satisfaire une condition de
félicité qui est de préserver certaines constantes. C'est ainsi qu'est
formulée dans l'\textsc{eme} une propriété fondamentale des chaines de
référence objective qui permet à l'information de circuler le long de
la chaine dans les deux sens. Pour les chaines démonstratives, la
condition de félicité est le transport de l'évidence: chaque saut est
lui-même un cheminement qui doit permettre à l'évidence de se propager
depuis l'énoncé initial source du saut jusqu'à l'énoncé final but du
saut.

Enfin, une dernière analogie est relative au caractère lointain,
c'est-à-dire éloigné des affaires communes, des entités amenées à
l'existence selon l'un ou l'autre des deux modes. En effet, il ne faut
rien de moins qu'une chaine --~de référence objective ou
démonstrative~-- et donc toute une série de médiations, d'instruments
et/ou de raisonnements, pour réussir à les instaurer.

Une source notable de dissemblance est relative au moteur de la
conviction obtenue à l'issue du parcours de ces chaines. Pour une
chaine de référence objective, on peut discerner deux principaux
facteurs de conviction: le premier est la confiance dans le protocole
suivi tout au long de la chaine pour contrôler le maintien des
constantes, le second est la simplicité perceptive de la dernière
inscription. Autrement dit, le premier facteur est de nature globale
puisqu'il concerne l'ensemble de la chaine, et il correspond à la
satisfaction de la condition de félicité correspondante; le second est
local puisqu'il est en rapport avec la dernière inscription
seulement. Dans le cas d'une chaine démonstrative, cette conviction
réside dans l'expérience de l'évidence du dernier énoncé, qui n'est
pas du même ordre que l'expérience de la simplicité perceptive de cet
énoncé final. Cette dernière inscription devient évidente pour un
mathématicien non pas du fait qu'il la trouve simple à appréhender
mais parce que cet acteur expérimente qu'un certain cheminement de
l'évidence aboutit à cette inscription. Ainsi, le facteur de
conviction d'une chaine démonstrative est avant tout global et il se
confond avec la satisfaction de la condition de félicité
correspondante. Néanmoins, on peut raisonnablement argüer que
l'expérience que constitue le parcours d'une chaine démonstrative
apporte une certaine familiarité avec sa dernière inscription qui, en
conséquence, acquiert une forme de simplicité perceptive\footnote{Pour
  prendre un exemple scolaire, l'inscription qu'est l'énoncé du
  théorème de Pythagore est perçue comme terriblement abstraite et
  complexe pour l'élève qui la découvre dans son cursus; mais si par
  la suite, cet élève réussit à apprivoiser une démonstration de ce
  résultat, cette inscription devient à ses yeux plus concrète et
  familière, bref, la perception qu'il en a se modifie dans le sens
  d'une plus grande simplicité. Pour invoquer James, une expérience
  entretient toujours des relations (qui sont elles-mêmes l'objet
  d'expériences) avec un contexte constitué d'autres
  expériences. Ainsi, la perception d'une inscription n'est pas une
  expérience isolée et elle mobilise immédiatement d'autres
  expériences qui lui donnent une partie de ses caractéristiques.}. Ce
n'est donc peut-être pas cette propriété de simplicité perceptive de
la dernière inscription qui constitue la différence la plus notable
entre les deux types de chaines. C'est au niveau des conditions de
félicité que la dissemblance est frappante. En effet, dans le cas
d'une chaine démonstrative, cette condition se confond avec
l'expérience de l'évidence de la dernière inscription. Cette évidence
signifie que l'acteur concerné voit littéralement, perçoit
directement, c'est-à-dire sans recours à un oracle, la validité du
résultat final. Cette signification de l'évidence est entièrement
dépendante de l'aspect diagrammatique de la pensée mathématique. En
effet, cette vision directe du résultat final n'est rien d'autre que
l'expérience perceptive de diagrammes convenablement (et souvent
laborieusement) modifiés. Au contraire, une chaine de référence
objective inclut généralement des transformations --~utilisation
cruciale d'appareillages technologiques, généralisation inductive du
résultat d'expériences singulières~-- qui ne peuvent pas produire la
même forme d'évidence, même si la conviction qui en découle peut elle aussi être
très forte. En conclusion, même si elles sont voisines, les
conditions de félicité propres aux deux types de chaines ne peuvent
être confondues puisque le transport de l'évidence d'une chaine
démonstrative n'est fondamentalement pas de même nature que le
maintien de constantes d'une chaine de référence objective.

Le lecteur attentif mais peu familier du point de vue de
l'\textsc{eme} peut s'attendre à ce que l'on souligne une autre
différence «évidente» concernant le rapport au monde matériel des
entités amenées à l'existence dans chacun de ces modes, puisqu'il peut
être clair pour lui que les êtres de la référence objective se
caractérisent par une présence active dans un secteur de la réalité
matérielle, alors que les êtres mathématiques se cantonnent dans des
mondes symboliques et abstraits sans ancrage matériel. Il se trouve
que cette évidence est remise en cause par l'\textsc{eme}, du fait
même que la notion de monde matériel ou de matière y est analysée
comme une confusion --~un amalgame entre deux modes d'existence~--
d'où découle la pauvreté de l'ontologie des modernes
(voir~\cite[chapitre~4]{EME}). Cette notion de matière au sens des
modernes est profondément transformée dans l'\textsc{eme}; en
particulier, elle acquiert un statut multimodal: chaque mode
d'existence possède sa propre notion de matière qui désigne l'ensemble
des entités dont dépendent les êtres de ce mode. L'effondrement dans
le cadre de l'\textsc{eme} de la notion unimodale de monde matériel
fait que l'attente supposée du lecteur perd toute pertinence et ne
peut qu'être abandonnée.
    
En guise de conclusion toute provisoire, le quasi-mode d'existence des
entités mathématiques décrit dans la partie~\ref{quasimode} et le mode
de la référence objective présentent un nombre significatif
d'analogies qui rendent difficiles une séparation pure et simple: ces
modes ont, pour reprendre une notion chère au second Wittgenstein,
comme un air de famille. Cependant, ces deux mêmes modes exhibent
aussi au moins une forme de différence assez nette qui empêche
l'identification des deux ou même la fusion de l'un dans l'autre. Or,
il semble inévitable de pouvoir disposer de ces deux modes afin de
couvrir le maximum du champ scientifique. D'ailleurs, c'est le plus
souvent un mixte des deux qui semble être à l'œuvre lorsque la
démarche considérée présente à la fois un développement
mathématico-déductif et un dispositif expérimental. Un exemple récent
et spectaculaire en est l'extraordinaire quête étalée sur près
d'un siècle qui vient de confirmer l'existence des trous noirs
\emph{via} la perception des ondes gravitationnelles\footnote{Voir la
  notion de «treuil ontologique» introduite par Étienne Klein
  \citep{klein08} pour évoquer le rôle des mathématiques dans
  l'ontologie de la physique. Le risque de ranimer l'antienne des
  mathématiques comme langage de la nature est implicitement présent,
  mais l'ancrage dans le cadre des modes d'existence de l'\textsc{eme}
  ajouté au quasi-mode présenté dans la présente contribution devrait
  constituer un puissant antidote.}. La démarche mythique de Thalès
déterminant la hauteur d'une pyramide par l'intermédiaire du théorème
qui porte aujourd'hui son
nom est tout aussi exemplaire. À la réflexion,
cette idée d'associer les deux modes semble inévitable pour rendre
compte de l'activité scientifique même dans le cas des sciences
expérimentales ou de terrain. En effet, puisque les inscriptions qui
constituent les étapes des chaines de référence font fortement
intervenir des entités mathématiques, le moindre traitement
mathématique de ces inscriptions --~par exemple un simple calcul
numérique~-- en toute rigueur relève du quasi-mode d'existence des
êtres mathématiques.  Pour tenir compte de ces contraintes, une
solution est de concevoir le mode [\textsc{ref}] sous la forme d'un
mode élargi contenant non seulement le mode de la référence objective
restreint décrit dans la partie~\ref{PrecisREF}, mais aussi le
quasi-mode des êtres mathématiques\footnote{Remarquons que les mathématiques tributaires du calcul sur machine évoquées dans la partie~\ref{DeuxRem} ont un mode
  d'existence différent des mathématiques habituelles mais relevant
  clairement de ce mode [\textsc{ref}] élargi.}. Cela revient à penser
[\textsc{ref}] comme un mode composé\footnote{À propos de cette idée
  de mode composé, il y a un précédent dans l'\textsc{eme} avec le cas
  du mode \emph{reproduction} [\textsc{rep}] qui se scinde lui-même en
  deux notions distinctes --~deux quasi-modes?~-- les \emph{inertes}
  et les \emph{vivants}.}.  Si cette hypothèse est retenue, il restera
à donner une description cohérente et unifiée de cette acception
étendue du mode [\textsc{ref}]. En attendant, il est facile
d'esquisser une description schématique de ce mode: ses
\emph{trajectoires} sont des chaines d'inscriptions obtenues chacune à
partir de la précédente par une transformation; son \emph{hiatus} est
la dissemblance, l'éloignement des inscriptions successives; sa
\emph{condition de félicité} est la circulation de la référence
objective ou de l'évidence le long de la chaine.
 
\begin{quote}
  \emph{Guy Wallet remercie Isabelle Stengers pour ses remarques et
    encouragements relativement à une version préliminaire de ce
    travail. Il remercie aussi Bruno Latour pour l'intérêt appuyé
    qu'il a manifesté pour une seconde version du même travail, et
    pour une discussion qu'il a organisée dans son laboratoire à ce
    propos. Enfin, les deux auteurs remercient Fabien Ferri pour les
    avoir mis sur la piste des thèses de Peirce commentés dans les
    écrits \citep{Chauvire1987,Chauvire2008} de Christiane Chauviré.}
\end{quote}

\noindent\textbf{Post-scriptum. } Tout en reconnaissant l'originalité et l'intérêt de l'analyse philosophique de la notion de démonstration développée dans le présent texte, un relecteur a fait part de son regret que le quasi-mode d'existence des êtres mathématiques qui en découle se cantonne à la seule «preuve d'énoncés déjà donnés» en laissant dans l'ombre la «généalogie des concepts». De manière générale, il oppose les «gestes à la Châtelet», comme l'invention en Grèce de la ligne sans épaisseur, par lesquels se transforment les paysages conceptuels, à la «monomanie», liée au logicisme, qui réduit les mathématiques à un amoncellement de démonstrations. Nous remercions ce relecteur qui, par cette critique intéressante, nous amène à préciser notre point de vue sur l'importance de l'innovation conceptuelle en mathématiques.

Cette innovation se manifeste dans le travail créatif de la preuve. Un exemple emblématique en est donné par les manipulations et enrichissements de diagrammes à l'œuvre dans les constructions auxiliaires de la géométrie des Grecs. D'après l'analyse de Peirce sur la dimension diagrammatique des mathématiques, il s'agit là d'un caractère général et consubstantiel de cette science, valant bien au delà de la seule géométrie. À la suite des travaux de Netz et de Peirce, notre analyse fournit un cadre radicalement éloigné du logicisme et du formalisme pour penser ces actes créatifs et les insérer dans la pratique du travail démonstratif.

\sloppy L'autre forme d'inventivité conceptuelle est celle qui surgit en amont des preuves, par exemple dans une définition, dans une conjecture, dans le dévoilement d'une analogie entre deux domaines \emph{a priori} sans rapport, etc. Ces gestes contribuent à modifier les champs d'étude et la compréhension que nous en avons. Pour cela, ils doivent affecter l'activité mathématique par excellence depuis les Grecs qui est le travail d'élaboration des preuves de nouveaux jugements. En accord avec le pragmatisme de Peirce et de James, nous considérons que, \emph{in fine}, l'importance et le sens d'un tel geste se mesure à ces traces concrètes. Il n'est pas dans l'objectif et la philosophie du présent article d'en interroger le sens lors de l'apparition initiale de ce geste en tant que «stratagème allusif» \citep[p.~221]{chatelet94}, avant toutes ces conséquences.

L'ontologie des mathématiques que nous avons développée est soucieuse d'une part de la pratique quotidienne dans les lieux où se développent et se communiquent les mathématiques et d'autre part des conditions concrètes, cognitives ou matérielles, de cette pratique. Puisque cette pratique est polarisée par le travail concret sur les démonstrations, nous maintenons qu'il s'agit d'un bon niveau pour juger de l'ontologie de cette science, à condition toutefois de s'appuyer sur la conception empiriste de la démonstration développée dans la partie~\ref{NonForm}, qui est, faut-il le rappeler, sans rapport avec la caricature logiciste ou formaliste de la preuve.  Cette dernière étant la cible explicite de la critique formulée par notre relecteur, nous estimons que notre présent travail n'est pas pleinement concerné par elle, tout en reconnaissant l'intérêt de la démarche qui la porte.
 
\providecommand\urlprefix{}
\bibliographystyle{phsc3-fr}
\bibliography{BiblioEME}

\end{document}